\documentclass{amsart}
\usepackage{amssymb,stmaryrd}
\usepackage{amsfonts}
\usepackage{amstext}
\usepackage{algorithmic}
\usepackage{algorithm}
\usepackage{graphicx}
\usepackage{epstopdf}
\usepackage[all]{xy}
\usepackage{MnSymbol}

\numberwithin{equation}{section}
\newtheorem{theorem}{Theorem}[section]
\newtheorem{lemma}[theorem]{Lemma}
\newtheorem{proposition}[theorem]{Proposition}

\theoremstyle{definition}

\newtheorem{example}[theorem]{Example}

\theoremstyle{remark}

\newcommand{\R}{{\mathbb{R}}}

\newcommand{\C}{{\mathbb{C}}}
\newcommand{\Z}{{\mathbb{Z}}}

\newcommand{\CL}{{\mathcal{L}}}

\newcommand{\wedgeq}{{\wedge\kern-5pt\cdot}}

\newcommand{\tens}{\otimes}

\newcommand{\id}{{\rm id}}

\newcommand{\extd}{{\rm d}}
\newcommand{\del}{{\partial}}
\newcommand{\eps}{\epsilon}

\def\lcross{{>\!\!\!\triangleleft}}

\newcommand{\cj}{{\mathfrak{j}}}
\newcommand{\ci}{{\mathfrak{i}}}
\newcommand{\lo}{\llbracket}
\newcommand{\lc}{\rrbracket}
\newcommand{\Vol}{{\rm Vol}}

\begin{document}

\title[Quantum Koszul formula]{Quantum Koszul formula on quantum spacetime}
\keywords{noncommutative geometry, quantum Riemannian geometry, quantum gravity, codifferential, Hodge Laplacian, central extension, differential graded algebra, differential form}

\subjclass[2000]{Primary 81R50, 58B32, 83C57}

\author{Shahn Majid  \& Liam Williams}
\address{School of Mathematics, Mile End Rd, London E1 4NS, UK}

\email{s.majid@qmul.ac.uk, liam.williams@qmul.ac.uk}


\begin{abstract}Noncommutative or quantum Riemannian geometry has been proposed as an effective theory for aspects of quantum gravity. Here the metric is an invertible bimodule map $\Omega^1\tens_A\Omega^1\to A$ where $A$ is a possibly noncommutative or `quantum'  spacetime coordinate algebra and $(\Omega^1,\extd)$ is a specified bimodule of 1-forms or `differential calculus' over it. In this paper we explore the proposal of a `quantum Koszul formula' in \cite{Ma:rec} with initial data a  degree -2 bilinear map $\perp$  on the full exterior algebra $\Omega$ obeying the 4-term relations
\[ (-1)^{|\eta|} (\omega\eta)\perp\zeta+(\omega\perp\eta)\zeta=\omega\perp(\eta\zeta)+(-1)^{|\omega|+|\eta|}\omega(\eta\perp\zeta),\quad\forall\omega,\eta,\zeta\in\Omega\]
and a compatible degree -1 `codifferential' map $\delta$. These provide a quantum metric and interior product and a canonical bimodule connection $\nabla$ on all degrees. The theory is also more general than classically in that we do not assume symmetry of the metric nor that $\delta$ is obtained from the metric. We solve and interpret the $(\delta,\perp)$ data  on the bicrossproduct model quantum spacetime $[r,t]=\lambda r$ for its two standard choices of $\Omega$. For the $\alpha$-family calculus the construction includes the quantum Levi-Civita connection for a general quantum symmetric metric, while for the more standard $\beta=1$ calculus we find the quantum Levi-Civita connection for a quantum `metric' that in the classical limit is {\em antisymmetric}.
\end{abstract}
\maketitle 

\section{Introduction}

Noncommutative differential geometry (NCDG) has been proposed for some three decades now as a natural generalisation of classical differential geometry that does not assume that the coordinate algebra or their differentials commute. There are many motivations and applications, many of them still unexplored (eg to actual quantum systems) but one of them is now widely accepted as an important role, namely as an effective theory for quantum gravity effects expressed as quantising spacetime itself. Of historical interest here was \cite{Sny} in the 1940's, although this did not propose a closed spacetime algebra exactly but an embedding of it into something larger. Specific proposals relating to quantum gravity (the `Planck scale Hopf algebra') appeared in \cite{Ma:pla} where they led to one of the two main classes of quantum groups to emerge in the 1980s as well as to one of the first and most well-studied quantum spacetimes with quantum group symmetry, namely the Majid-Ruegg `bicrossproduct model' \cite{MaRue}. In 2D this is the coordinate algebra  $[r,t]=\lambda r$ where $\lambda$ should be $\imath$ times the Planck scale of around $10^{-35}$m. In spite of many hundreds of papers on this quantum spacetime, it continues to be useful  as a testbed for new ideas in noncommutative geometry and continues to surprise. In particular, it was shown recently in \cite{BegMa:cqg} that the standard differential calculus on this algebra, namely \begin{equation}\label{beta1calc} [r,\extd t]=\lambda\extd r,  \quad [t,\extd t]=\lambda\extd t,\quad [r,\extd r]=0,\quad [t,\extd r]=0\end{equation}
admits only a 1-parameter form of quantum metrics which classical $\lambda\to 0$ limit, namely
\[  \extd r^2+ B v^2;\quad  v=r \extd t - t \extd r \]
which is that of either for $B>0$ an expanding universe with an initial big bang singularity {\em or} for $B<0$ a gravitational source so strong that even light eventually gets pulled back in and with a curvature singularity at $r=0$. The calculus here is the $\beta\ne 1$ point of a family of calculi with similar features. Then in \cite{MaTao} it was shown that the other $\alpha$ family choice of calculus similarly admits a unique form of quantum metric which is either de Sitter or anti-de Sitter space depending on the sign of a parameter. Up to a change of variables we can again take $\alpha=1$, then
\begin{equation}\label{alpha1calc} [t, \extd r]=-\lambda\extd r, \quad [t,\extd t]=\lambda\extd t \end{equation}
is the calculus, and the quantum metric has classical limit
\[ r^{-2}\extd r^2+ 2 a\extd r\extd t+b r^{2}\extd t^2. \]
with $a^2>b$. The classical geometry here depends on the sign of $b$. In both cases we see that a particular classical (pseudo)Riemannian geometry {\em emerges} as being forced out of nothing but the choice of algebra and its differential structure, showing that the `quantum spacetime hypothesis' has implications for classical GR. These constraints on classical geometry emerging from noncommutative algebra were analysed in general at the semiclassical level, as a new theory of Poisson-Riemannian geometry, in \cite{BegMa:poi}. Moreover, in both cases the full quantum geometry is constructed in \cite{BegMa:cqg,MaTao} in the sense of a quantum-Levi Civita (or quantum torsion free quantum metric compatible) connection in the bimodule formalism of quantum Riemannian geometry in that has its roots in \cite{DV1,DVM,Mou,BegMa:com}. 

In spite of these successes, the general formalism of `quantum Riemannian geometry' in both the bimodule connection approach and an earlier quantum group frame bundle approach\cite{Ma:rie} has until now lacked a general construction for the  quantum-Levi-Civita connection, which has to be solved for on a case by case basis. Recently in \cite{Ma:rec}, however, one of the present authors introduced a radically new point of view on both classical and quantum Riemannian geometry as emerging from a choice of codifferential $\delta$ (not the other way around as would be more usual) along with a new formula\cite{Ma:rec}
\begin{equation}\label{kos}  \nabla^{LC}_\omega\eta={1\over 2}\left(L_\delta(\omega,\eta)+\CL_\omega\eta+ (\extd\omega)\perp\eta\right);\quad L_\delta(\omega,\eta)=\delta(\omega\eta)-(\delta\omega)\eta+\omega\delta\eta\end{equation}
for the classical Levi-Civita connection. Here we view a 1-form $\omega\in\Omega^1(M)$ as a vector field via the metric and $(\ )\perp\eta$ is interior product by the vector field similarly corresponding to $\eta\in\Omega^1(M)$. The Lie derivative is also given by such an interior product $\omega\perp$ and $\extd$. The work also  led to a new property\cite{Ma:rec}
\[  \delta(\omega\eta\zeta)=(\delta(\omega\eta))\zeta+(-1)^{|\omega|}\omega\delta(\eta\zeta)+(-1)^{(|\omega|-1)|\eta|}\eta\delta(\omega\zeta)\]
\[\qquad\qquad- (\delta\omega)\eta\zeta-(-1)^{|\omega|}\omega(\delta\eta)\zeta-(-1)^{|\omega|+|\eta|}\omega\eta\delta\zeta\]
for the classical codifferential acting on  $\omega,\eta,\zeta\in \Omega(M)$, the exterior algebra on the manifold. This says, remarkably, that $(\Omega(M),\delta)$ makes any Riemannian manifold into a Batalin-Vilkovisky algebra. From our new starting point we can go further and axiomatise $\delta$ as a degree -1 map obeying certain axioms and if this is of `classical type' (notably $\delta^2$ is tensorial, for example  zero) then the connection defined as above will  necessarily be the Levi-Civita one for an inverse metric $(\ ,\ )$ induced by $\delta$ according to the formula
\[ \delta(f\omega)=f\delta\omega+(\extd f,\omega),\]
for all $f\in C^\infty(M),\omega\in \Omega^1(M)$, see \cite[Thm 3.18]{Ma:rec}. Another feature of this new approach to classical Riemannian geometry is that it works well with forms of all degree. Thus the above formula for $\nabla^{LC}$ works for $\eta$ of all degrees provided we extend $\perp$ to all degrees by the formula\cite{Ma:rec}
\[ (\omega_1\cdots\omega_m)\perp(\eta_1\cdots\eta_n)=\sum_{i,j}(-1)^{i+j}(\omega_i,\eta_j)\omega_1\cdots\widehat{\omega_i}\cdots\omega_m\eta_1\cdots \widehat{\eta_j}\cdots\eta_n,\quad\omega_i,\eta_j\in \Omega^1(M),\]
where we leave out the hatted ones. If $\omega$ has degree 1 then $\omega\perp(\ )$ is interior product as used in the Lie derivative  in the formula for $\nabla^{LC}$. Classically, $\perp$ is not more data than the metric, it merely extends it as a bi-interior product, and our Koszul formula is equivalent in this case to the usual Koszul or Levi-Civita formula but in a novel differential form language that depends also on constructing the associated Hodge codifferential $\delta$ compatibly with the metric\cite{Ma:rec}. On the other hand, even when $A=C^\infty(M)$, we are not limited to this choice as we could let $\perp$ be nonsymmetric and still define the inverse metric as the symmetrisation of $\perp$ in degree 1 in the construction of \cite{Ma:rec}, and we are also not limited to the standard `classical type' $\delta$ (we look at this slightly more general but still  classical construction in Section~2.3). 

It was also pointed out but not the main topic of \cite{Ma:rec} that this differential-Koszul formula can be applied when our algebra of coordinates is a noncommutative algebra $A$ to begin with, and $(\Omega(A),\extd)$ is a quantum differential calculus. We still need a map $\perp$ which we axiomatise as a degree -2 `product' $\perp$ on  $\Omega(A)$ obeying\cite{Ma:rec}  
\[ (-1)^{|\eta|} (\omega\eta)\perp\zeta+(\omega\perp\eta)\zeta=\omega\perp(\eta\zeta)+(-1)^{|\omega|+|\eta|}\omega(\eta\perp\zeta),\quad\forall\omega,\eta,\zeta\in\Omega(A)\]
(which we call the `4 term relation') together with a degree $-1$ map $\delta:\Omega(A)\to \Omega(A)$ the `quantum codifferential' such that 
\[  \delta(a\omega)=a\delta\omega+\extd a\perp\omega,\quad \delta(\omega a)=(\delta\omega)a+\omega\perp^R \extd a\]
for a $\perp$ as above and another similar bimodule map $\perp^R$. These formulae determine $\perp,\perp^R$ on degree 1 if we take $\delta$ as a starting point. From this data, it is shown that one can construct a  quantum bimodule connection $\nabla$ from $\perp$ by the same formula (\ref{kos}), with quantum (inverse) metric $(\ ,\ )={1\over 2}(\perp+\perp^R)$ when restricted to 1-forms. Classically $\perp^R=\perp\circ{\rm flip}$ so this would be symmetric. We also obtain a `quantum interior product' $\cj$ by allowing higher degree forms in the first argument, which is also something lacking in noncommutative geometry.   The quantum connection (\ref{kos}) now is not necessarily torsion free and quantum metric compatible or `quantum Levi-Civita' (or QLC) in the sense of \cite{BegMa:cqg} but in so far as we make choices that deform the classical theory, the connection will deform the classical $\nabla^{LC}$.  Moreover, the construction has its own interest which does still apply in the quantum or noncommutative case, and which we explain next. This makes these quantum bimodule connections natural and of interest in their own right even if they do not  necessarily obey exactly the previously proposed axioms of a QLC, in which case deviation from the latter would now be viewed as a source of new effects. 

Specifically, the $(\delta,\perp)$ construction arises in \cite{Ma:rec} much more deeply from nothing but the axioms of a noncommutative differential calculus (basically, the Leibniz rule) and a central extension problem. Thus, in the classical case, one can look for 
\[ \Omega_{\theta'}\hookrightarrow\tilde\Omega\twoheadrightarrow \Omega(M)\]
as a sequence of differential graded algebras where  we extend the classical exterior algebra to a quantum one $\tilde\Omega$ by adjoining  a graded-commuting $\theta'$ with $\extd\theta'=0$ and $\theta'^2=0$. Such an extension is called `cleft' if $\Omega\cong \Omega(M)\oplus\theta'\Omega(M)$ as a left $C^\infty(M)$-module and `flat' if it is equivalent to a cleft extension with $\extd$ undeformed. It was shown in \cite{Ma:rec} that cleft central extensions are in 1-1 correspondence with certain 2-cocycle data $([[\ ,\ ]],\Delta)$ that can be interpreted as including a  possibly degenerate (pseudo) Riemannian metric $(\ ,\ )$ as part of an interior product map $\cj$, a connection $\nabla$ and a Laplacian. In the flat case $\Delta=\extd\delta+\delta\extd$ for some codifferential $\delta$ and $\nabla$ is the Levi-Civita connection given by (\ref{kos}), on all degrees. This gives a mechanism by which the structures of classical GR could emerge out of the algebraic structure of quantum spacetime if its quantum  differential calculus approaches a central extension as we approach the classical limit. One reason why this could typically be the case is what has been called the `quantum anomaly for differentials' in the quantum group literature: often there is not a suitably covariant differential calculus within deformation theory (due to the the lack of a flat covariant Poisson connection from the point of view of Poisson-Riemannian geometry\cite{BegMa:poi}) and one must either live with a nonassociative differential calculus or absorb the anomaly by having a higher dimension\cite{BegMa:sem}. The same extension theory as above applies when we replace $\Omega(M)$ by some quantum $\Omega(A)$ and a flat cleft extension of that leads to both $\delta$ and a cocycle $([[\ ,\ ]],\Delta)$ which is shown in \cite{Ma:rec} to provide a bimodule connection when the first argument of the bracket is in $\Omega^1(A)$ as well as an interior product $\cj$ when the second argument is degree zero. Details are in Section~2.1.

Thus we have a deeper point of view on how the familiar structures of GR could arise purely out of noncommutative differential algebra, as well as a practical route to quantise them. In the present paper we will explore these new ideas in the context of the bicrossproduct quantum spacetime $[r,t]=\lambda r$ with its two choices (\ref{alpha1calc}) and (\ref{beta1calc}) of differential calculi. In both cases one has a basis $\{e_i\}$ of central 1-forms (that commute with all functions) and an inverse quantum metric $g^{ij}=(e_i, e_j)$ as any $2\times 2$ constant matrix of coefficients (we do not impose quantum symmetry or `reality' conditions as in \cite{BegMa:cqg,MaTao} so do we not have a unique form of metric).  We also could have any constant matrix for the coefficients for the interior product $\cj_{e_i}(\Vol)=v^{ij}e_j$ where $\Vol=e_1e_2$ is the central top form.  In Section~3 to solve the 4-term relations with differentials (\ref{alpha1calc}) to find that $b_{ij}=e_i\perp e_j$ is any $2\times 2$ matrix with constant entries. We then take a general form of $\delta$,  and apply the Koszul formula to construct a quantum bimodule connection $\nabla$, metric $g$ and interior product $\cj$. Remarkably, we find that some of the conditions for $\delta^2$ to a left-module map or `left-tensorial' precisely characterise a class of quantum connections with classical limit as $\lambda\to 0$, see Theorem~\ref{thmalphalc}. Among this class and for generic $b_{ij}$, we find:
\begin{enumerate}
\item[(i)] The interior product is  $v^{ij}=g^{ki}\eps_{kj}$ as classically, where $\eps_{12}=1$ is the antisymmetric tensor;
\item[(ii)] $\nabla$ is then QLC, i.e. torsion free and quantum metric compatible, if and only if $\delta^2$ is a `strongly tensorial' in the sense of a bimodule map;
\item[(iii)] The metric needed for this  is  $g^{ij}=(b_{ij}+b_{ji})/2$, as classically.
\item[(iv)] The $\delta$ needed form a two parameter space of constant $a_i, b_i$ with $b_i$ determined, including the case where $\delta^2=0$ as classically.
\end{enumerate}
The quantum Koszul formula in this case works as expected. It not only gives the previously known connnection\cite{MaTao} but adds the interior product and `explains' why the metric that emerges is symmetric rather than this being assumed as in \cite{MaTao},  namely in order to be compatible with the connection induced by the quantum central extension data.

In Section~4 we similarly solve the 4-term relation for the same quantum spacetime and its `standard'  differentials (\ref{beta1calc}). This time we  find a unique form of $\perp$ namely $b_{ij}=b\begin{pmatrix}0 &1\cr -1 & -\lambda\end{pmatrix}$, which we see has an unexpected antisymmetric  form in the classical limit as $\lambda\to 0$. We again find that $\delta^2$ a left-module map ensures a classical limit for the connection given by the quantum Koszul formula, see Theorem~\ref{thmbetalc} and requiring $\delta^2$ to be a bimodule map makes the  connection weak QLC. With a small further condition on the metric it becomes QLC, see Example~\ref{betaglc},
\begin{enumerate}
\item[(i)] The interior product is $v^{ij}=\begin{pmatrix}-g^{21} &g^{11}\cr g^{22} & g^{12}\end{pmatrix}=g^{ki}\eps_{kj}+O(\lambda)$;
\item[(ii)] $\nabla$ is  QLC;
 \item[(iii)] The metric needed has the form $g^{ij}=g^{12}\begin{pmatrix}0 &1\cr -1 & -{\lambda\over 2}\end{pmatrix}=g^{12}\eps_{ij}+O(\lambda)$;
\item[(iv)] The $\delta$ needed has an order $1/\lambda$ singularity as $\lambda \to 0$, is uniquely determined up to a constant of integration and has $\delta^2=0$. 
\end{enumerate}

In both cases we  can land on any freely chosen $g^{ij}$ by choice of $(\delta,\perp)$ and we can further choose $\delta^2$ a left module map, which ensures classical limits and that $v^{ij}$ is built from $g^{ij}$, but without further restrictions $\nabla$  need not be torsion free or quantum metric compatible  or even a weaker `cotorsion free' version of the latter\cite{Ma:rie} which is common in noncommutative Riemannian geometry.  In both Sections~3,4 we provide a rather fuller analysis of the properties resulting from different assumptions on $\delta$, including results motivated from a general feature of connections coming from central extensions of classical type in \cite[Prop 3.16]{Ma:rec}  whereby the torsion and metric compatibility are linearly related, but now in our quantum examples. The quantum symmetric metric in \cite{BegMa:cqg} is then covered in Example~\ref{betabh} for which we obtain in the limit a particular classical connection which is not the Levi-Civita one, and quantise it. The paper ends with some concluding remarks.

\section{preliminaries}

Throughout the paper, a differential graded algebras or DGA over an algebra $A$ means a graded algebra $\Omega=\oplus_n\Omega^n$ with $\Omega^0=A$ and $\extd:\Omega^n\to \Omega^{n+1}$ for all available degrees with $\extd^2=0$ and $\extd$ obeying the graded Leibniz rule. We will say that a DGA is {\em standard} (or an `exterior algebra') if generated by $A,\extd$. 

By {\em quantum (inverse) metric} we mean that $\Omega^1$ is equipped with a bimodule map $(\ ,\ ):\Omega^1\tens_A\Omega^1\to A$ and normally we will assume this is invertible so there is an actual element $g=g^1\tens g^1\in \Omega^1\tens\Omega^1$ (sum understood) inverse to it in the sense $(\omega, g^1)g^2=\omega=g^1(g^2,\omega)$ for all $\omega\in\Omega^1$. As shown in \cite{BegMa:cqg} this will entail that $g$ is central. However, in the present paper $(\ ,\ )$ appears to play a more important role and we may allow it to be degenerate. In our construction it appears as the degree 1 case of a more general quantum interior product $\cj:\Omega^n\tens_A\Omega^1\to\Omega^{n-1}$ and we sometimes denote it by $\cj$ for this reason. 

By a {\em form-covariant derivative} we mean  $\nabla_\omega:E\to E$ where $E$ is a left $A$-module and 
\[\nabla_\omega(a.e )= \nabla_{\omega a}(e)+ (\omega,\extd a)e,\quad \nabla_{a\omega}=a\nabla_\omega\]
which is based on the usual axioms in noncommutative geometry for a left connection but evaluated against a 1-form via $(\ ,\ )$. We have a `bimodule covariant derivative' if $E$ is a bimodule and there is a bimodule map $\sigma:\Omega^1\tens_A E\tens_A\Omega^1\to E$ such that
\[ \nabla_\omega(e.a)=(\nabla_\omega e).a+ \sigma_\omega(e\tens\extd a)\]
which is evaluation against $\cj$ of the usual notion \cite{DV1,DVM,Mou,BegMa:com} of a bimodule covariant derivative. Moreover, if $\cj$ is invertible with inverse $g=g^1\tens g^2$ then $\nabla_\omega$  on $E,F$ has a tensor product 
\begin{equation}\label{covtens} \nabla_\omega(e \tens_A f)=\nabla_\omega e\tens_A f+ \sigma_\omega(e \tens_A g^1 )\tens_A\nabla_{g^2  }f,\quad\forall e\tens f\in E\tens_A F.\end{equation}

Setting $E=F=\Omega^1$ we can ask that $\nabla_\omega(g)=0$ which is the notion of a bimodule connection on $\Omega^1$ being {\em metric compatible}. In terms of $(\ ,\ )$, if $\sigma$ is invertible, then metric compatibility is equivalent to 
\[ (\id\tens (\ ,\ ) )\nabla(\omega\tens\eta)=\extd (\omega,\eta),\quad\forall \omega,\eta\in \Omega^1.\]
Also when $E=\Omega^1$ and $(\ ,\ )$ is invertible, one has the notion of {\em torsion}, 
\[ T_\nabla:\Omega^1\to\Omega^2,\quad T_\nabla=\extd - g^1\nabla_{g^2}\]
and the notion of {\em cotorsion}
\[ coT_\nabla=(\extd\tens\id -\id\wedge\nabla)g=\extd g^1\tens g^2-g^{1'}\nabla_{g^{2'}}g^2\in \Omega^2\tens_A\Omega^1\]
for any connection on $\Omega^1$, where the primes denote a second copy of $g$. By definition a connection is {\em quantum Levi-Civita} if it is torsion free and metric compatible. It is {\em weak quantum Levi-Civitia} if it is torsion and cotorsion free (often in noncommutative geometry this weaker property is all we have).  A connection has curvature defined by
\[ R_\nabla=(\extd\tens\id-\id\wedge\nabla)\nabla: \Omega^1\to \Omega^2\tens_A\Omega^1\]
which can also be converted in terms of $\nabla_\omega$.  Apart from translating to form-derivatives, these  are all established notions of a constructive approach to noncommutative geometry, see \cite{DVM,Mou,BegMa:com,BegMa:cqg,BegMa:poi} and references therein.

\subsection{Central extensions} \label{SecExt}

The notion of a central extension $\tilde \Omega(A)$ of a DGA $\Omega(A)$  was introduced in \cite{Ma:rec} as an extension in degree 1 by the algebra $\Omega_{\theta'}=k[\theta']/\langle\theta'^2\rangle$ viewed as a trivial DGA with $\theta'$ of degree 1 and $\mathrm{d}\theta'=0$. More precisely,  $$\tilde{\Omega}(A)=\Omega_{\theta'}\otimes\Omega(A)$$ as a vector space and $$0\rightarrow\Omega_{\theta'}\hookrightarrow\tilde{\Omega}(A)\twoheadrightarrow\Omega(A)\rightarrow 0$$ as maps of DGA's, where the outer maps come from the canonical inclusion in the tensor product and by setting $\theta'=0$. We also require that $\theta'$ is graded-central, $$\theta'\omega=(-1)^{|\omega|}\omega\theta'$$
in $\tilde{\Omega}(A)$. A morphism of extensions $\Phi:\tilde{\Omega}(A)\rightarrow \tilde{\Omega}'(A)$ means a map of DGA's such that $$\Phi(\theta')=\theta',\quad \Phi(\omega)=\omega-\frac{1}{2}\theta'\delta(\omega)$$ By a (left) cleft extension we mean a central extension where the canonical linear inclusion of $\Omega(A)$ coming from the tensor product form is a left A-module map. And by a flat extension we mean one which is equivalent to one where $\extd$ is not deformed. 

It is shown in \cite{Ma:rec} that an extension must necessarily have the form $$\omega\cdot\eta=\omega\eta-\frac{1}{2}\theta'\llbracket\omega,\eta\rrbracket, \quad \mathrm{d}.\omega=\mathrm{d}-\frac{1}{2}\theta'\Delta\omega, \quad \omega,\eta\in\Omega(A)$$ for a bilinear map $\llbracket\; ,\; \rrbracket$ of degree -1 and a linear map $\Delta$ of degree 0, forming a `cocycle' in the sense\cite{Ma:rec}
\begin{equation}\label{c1} \lo\omega\eta,\zeta\lc+\lo\omega,\eta\lc\zeta=\lo\omega,\eta\zeta\lc+(-1)^{|\omega|}\omega\lo\eta,\zeta\lc.\end{equation}
\begin{equation}\label{c2} L_\Delta(\omega,\eta)=\extd\lo\omega,\eta\lc+\lo\extd\omega,\eta\lc+(-1)^{|\omega|}\lo\omega,\extd\eta\lc
\end{equation}
for all $\omega,\eta,\zeta\in \Omega(A)$, and $[\Delta,\extd]=0$. 
The extension is cleft precisely when $\llbracket a,\; \rrbracket=0$ for all $a\in A$ and is flat precisely when $\Delta=\extd\delta+\delta\extd$ for some degree -1 map $\delta$ \cite{Ma:rec}. 

A cleft extension $(\Delta,\lo\ ,\ \lc)$ on a standard DGA $\Omega(A)$ is $n$-{\em regular} if\cite{Ma:rec}
\[ \cj_\omega(a\extd b)={1\over 2}\lo \omega a,b\lc,\quad \forall \omega\in \Omega,\ a,b\in A\]
is  a well-defined degree -1 map $\cj:\Omega^i \tens_A\Omega^1\to \Omega^{i-1}$ for $i\le n$. We say that the cleft extension is regular if it is regular for all degrees. We refer to $\cj$ as `interior product' and its restriction $\Omega^1\tens_A\Omega^1\to A$  to degree 1 will be the `inverse metric' $(\ ,\ )$ of the geometry induced by a central extension.  It is shown in \cite[Prop~3.6]{Ma:rec} that if $(\Delta,\lo\ ,\ \lc)$ is a regular cleft extension on a standard DGA $\Omega(A)$ then $\cj$ is a bimodule map and
\[ \nabla_\omega\eta={1\over 2} \lo \omega,\eta\lc,\quad\forall \omega\in \Omega^1,\quad \eta\in\Omega\]
is a bimodule covariant derivative on $\Omega$ with respect to $\cj$. Here
\[ \sigma:\Omega^1\tens_A\Omega\tens_A\Omega^1\to \Omega,\quad \sigma_\omega(\eta\tens_A\zeta)=\cj_{\omega\eta}(\zeta)+\omega\cj_\eta(\zeta),\quad\forall \omega,\zeta\in \Omega^1,\ \eta\in \Omega.\] 
For $\nabla$ on $\Omega^1$ to be torsion free in the case of a standard calculus needs
\begin{equation}\label{torsionfreeg} g^1[[g^2,\extd a]]=0\end{equation}
and to be metric compatible, given the form of $\sigma$, needs
\begin{equation}\label{metriccompatg}(\nabla_\omega(g^1b_j)-(\nabla_\omega g^1)b_j)\tens_A\nabla_{g^j} g^2+\nabla_\omega g^1\tens_A g^2=0 \end{equation}
where $g=\extd b^j\tens_A g^j$ which we can write in terms of $[[\ ,\ ]]$. The weaker cotorsion free condition becomes
\begin{equation}\label{cotorsiong}(\extd b^j)g^1\tens \nabla_{g^2} g^j=0.\end{equation}
Note also that just as $\cj$ is not necessarily a derivation, we do not necessarily have compatibility of the connection on higher forms with the wedge product, i.e. the {\em braided Leibniz condition}\cite{Ma:rec}
\[ \nabla_\omega(\eta\zeta)=(\nabla_\omega \eta)\zeta+\sigma_\omega(\eta\tens g^1)\nabla_{g^2}\zeta\]
which using the cocycle condition would come down to
\begin{equation}\label{wedgecompatg} j_{\omega\eta}(g^1)[[g^2,\ ]]=[[\omega\eta,\ ]],\quad\forall \omega\in \Omega^1,\quad \eta\in\Omega.\end{equation}

\subsection{Quantum Koszul formula}\label{SecKos}

In \cite{Ma:rec} there is a construction for central extensions based on a bimodule map $\perp$ of degree -2 on $\Omega(A)$ which on degree 1 we view as an inverse metric $(\ ,\ ):\Omega^1\tens_A\Omega^1$ and which in general has to obey  $\perp a=a\perp=0$ for all $a\in A$ and  the {\em 4-term relation}
\[ (-1)^{|\eta|} (\omega\eta)\perp\zeta+(\omega\perp\eta)\zeta=\omega\perp(\eta\zeta)+(-1)^{|\omega|+|\eta|}\omega(\eta\perp\zeta),\quad\forall\omega,\eta,\zeta\in\Omega\]
which implies when one of the arguments is degree 0 that $\perp:\Omega^m\tens_A\Omega^n\to \Omega^{m+n-2}$ is a well defined bimodule map. 

\begin{theorem}\label{thm4term}\cite[Thm 3.12]{Ma:rec} If $\Omega(A)$ is a standard DGA, $\perp$ obeys the 4-term relation, and if $\delta$ is a degree -1 map such that 
\[ \delta(a\omega)-a\delta(\omega)=\extd a\perp\omega,\quad \delta(\omega a)-(\delta\omega)a=\omega\perp_R \extd a,\quad \forall a\in A,\ \omega\in\Omega\]
for some bimodule map $\perp_R:\Omega^n\tens_A\Omega^1\to \Omega^{n-1}$ (we say that $\delta$ is `regular'). Then we have a regular flat cleft extension with 
\[ \Delta=\extd\delta+\delta\extd,\quad \lo\omega,\eta\lc=L_\delta(\omega,\eta)+\omega\perp\extd\eta-(-1)^{|\omega|}\extd\omega\perp\eta-(-1)^{|\omega|}\extd(\omega\perp\eta),\quad\forall \omega,\eta\in \Omega.\]
\end{theorem}

By the above, this implies $\nabla$ a bimodule covariant derivative and candidate for a `quantum Levi-Civita'-like connection for the quantum metric $ g$. It also implies an interpretation of the interior product  $\cj_\omega(\extd a)={1\over 2}[[\omega,a]]$ as a `connection on degree 0' which from the above is
\[ \cj_\omega(\eta)={1\over 2}(\omega\perp\eta+\omega\perp_R\eta),\quad \eta\in\Omega^1,\ \omega\in\Omega\]
extend the inverse quantum metric $(\ ,\ )$ to all degrees in its first input. In the classical case it is shown in \cite{Ma:rec} that $\nabla$ indeed is torsion free and metric compatible with $g$ inverse to $(\ ,\ )=\perp=\perp_R$ if we take for $\delta$ the standard Riemannian codifferential, and then $\cj_\cdot(\eta)$ is indeed the interior product along the vector field corresponding via the metric to $\eta$. Note that the centrally extended noncommutative DGA $\tilde\Omega$ behind the theorem need not be standard.

It is shown in \cite{Ma:rec} that when a calculus is inner in the (purely `quantum') sense that there exists a 1-form $\theta\in \Omega^1$ such that  $\extd =[\theta,\ \}$ is the graded commutator, then 
\[ \delta=\theta\perp,\quad \perp_R=0\]
provides the required data for any solution of the 4-term relations and gives\cite{Ma:rec}
\[  \cj_\omega(\eta)={1\over 2}\omega\perp\eta,\quad \nabla_\omega=-{1\over 2}L_{\perp\theta}(\omega,\ ),\quad \sigma_\omega(\eta\tens\zeta)={1\over 2}\left((\omega\eta)\perp\zeta-(-1)^{|\omega|}\omega(\eta\perp\zeta)\right)\]
\[ \Delta=2\nabla_\theta-\theta^2\perp  \]
on forms of all degrees (on degree 0 this is $\Delta a=2\cj_\theta(\extd a)=\theta\perp\extd a=\delta\extd a$). One can check that $\nabla_\omega$ is evaluation by $\cj$ of $\nabla=\theta\tens\eta-\sigma(\eta\tens\theta)$. One can show that in general
\[ T_\nabla(\omega)=-\omega\theta-{1\over 2}g^1((g^2\omega)\perp\theta)=-2\omega\theta+{1\over 2}g^1\wedge (g^2\perp(\omega\theta))\]
\[ {\rm co}T_\nabla=2\theta g-{1\over 2}g^{1'}g^1\tens (g^2g^{2'})\perp\theta=3\theta g+{1\over 2}g^1g\perp(g^2\theta).\]
However, this is just one (far from classical) example $\delta$  and we may be more interested in prescribing $\cj$ to a given quantum metric and choosing $\delta$ as needed for this.

Finally, we remark that for a cocycle built in this way from data $(\delta,\perp)$ we actually have a further
extension $\widetilde\widetilde\Omega(A)\to\widetilde\Omega(A)\to\Omega(A)$ where we allow $\extd \theta'\ne 0$ namely with the new operations \cite[Prop~3.21]{Ma:rec}
\[ \omega\cdot\eta=\omega\eta+{\mu\over 2}(-1)^{|\omega|+|\eta|}\lo\omega,\eta\lc\theta'-{\mu\over 2}(-1)^{|\omega|}(\omega\perp\eta)\extd\theta' \]
\[ \extd_\cdot\omega=\extd \omega - {\mu\over 2}(-1)^{|\omega|}(\Delta\omega)\theta'+{\mu\over 2}(\delta\omega)\extd\theta',\quad\forall \omega,\eta\in \Omega(A),\]
and $\theta'^2=\theta'\extd\theta'=(\extd\theta')\theta'=\{\omega,\theta'\}=0$.

\subsection{Discrete nonommutative example}\label{SecZ2}

Although not our main topic, the theory applies to the commutative coordinate algebra  $A=\C(\Z_2\times\Z_2)$ with its direct product noncommutative differential calculus (each $\Z_2$ has a unique calculus, the universal one). Differential 1-forms on a discrete set can be identified as edges of a graph and this is the calculus on a square and we are solving for the noncommutative geometry of a square. 

The calculus has basis of translation invariant 1-forms $e_i$, $i=1,2$ with relations $e_i f= R_i (f)e_i$ where $R_i(f)$ is right translation in the $i$'th factor.  The exterior derivative on degree 0 is $\extd f=(\del_if)e_i$ where $\del_i=R_i-\id$. The exterior algebra is this model is defined in the usual way by $e_i^2=0$ and $e_1e_2+e_2e_1=0$, with top form $\Vol=e_1e_2$. For the map $\perp$ we are forced to take a diagonal form 
\[ e_i\perp e_j=\delta_{ij}a_i\]
since the bimodule relations require that $e_1\perp e_2 f=R_1R_2(f) e_1\perp e_2$ for all $f$ which since $e_1\perp e_2$ is an element of a commutative algebra is not possible unless it is zero. 

The 4-term relation on this DGA in degrees 1 on the diagonal case $e_i,e_i,e_i$ is
\[ a_i e_i=e_i a_i,\quad{\rm i.e.}\quad \del_i a_i=0\]
while if $i\ne j$ we have 
\[ -e_i^2\perp e_j+a_i e_j=\eps_{ij} e_i\perp\Vol,\quad  -\eps_{ij}\Vol\perp e_j=e_i\perp e_j^2+e_i a_j,\quad  -\eps_{ij}\Vol\perp e_i=\eps_{ji}e_i\perp\Vol \]
which means
\[ e_i\perp\Vol=a_i\eps_{ij}e_j,\quad \Vol\perp e_j=e_i\eps_{ij}a_j,\quad R_j(a_i)=-a_i.\]
If one of the forms is $\Vol$ then we have
\[ -\eps_{ij}\Vol\perp\Vol +\delta_{ij}a_j\Vol=e_i(e_j\perp\Vol),\quad -(\Vol\perp e_i)e_j=\Vol\perp\Vol \eps_{ij}-\Vol\delta_{ij}a_j\]
\[ (e_i\perp\Vol)e_j=-e_i(\Vol\perp e_j)\]
which when $i=j$ are all obeyed given the relations on the $a_i$. (One of these is $\Vol a_i = -a_i \Vol$.) When $i\ne j$ the first two are both equivalent to $\Vol\perp\Vol=0$ while the last is empty. We thus solve our 4-term relations with two constant parameters 
\[ a_1=(a,-a,a,-a),\quad a_2=(b,b,-b,-b),\quad a,b\in \C\]
where we list the values at the points 00,01,10,11 of $\Z_2\times\Z_2$. 

Because this is only a warm up, we will not do the full analysis of all
possible $\delta$ compatible with the above $\perp$, but merely give an example:

\begin{example} Up to an overall normalisation there is a unique 1-parameter form of quantum metric on $A=\C(\Z_2\times \Z_2)$ coming out
of the $\perp$ construction and quantum Koszul formula with $\delta=\theta\perp$, namely
\[ g=e_1{2\over a_1}\tens e_1+e_2{2\over a_2}\tens e_2\]
\[ \nabla e_1=-2\alpha e_2\tens e_2,\quad \nabla e_2=-2\alpha^{-1}e_1\tens e_1;\quad\alpha={a_1\over a_2}={a\over b}(1,-1,-1,1)\]
\[ \sigma(e_1\tens e_1)=e_1\tens e_1+2\alpha e_2\tens e_2,\quad \sigma(e_1\tens e_2)=e_2\tens e_1\]
\[  \sigma(e_2\tens e_1)=e_1\tens e_2,\quad \sigma(e_2\tens e_2)=e_2\tens e_2+2\alpha^{-1}e_1\tens e_1\]
which is invertible and not involutive. Here $\nabla e_i=\theta\tens e_i-\sigma(e_i\tens \theta)$ and is torsion-free and cotorsion-free
(or weak quantum Levi-Civita). It has curvature
\[ R_\nabla(e_1)=4\Vol\tens(\alpha e_2-e_1)=\alpha R_\nabla(e_2).\]
\end{example}
\proof The calculus here is inner with $\theta=e_1+e_2$ which gives in our case 
\[ \delta (e_i)=a_i,\quad \delta(\Vol)=a_1 e_2-a_2 e_1,\quad \nabla_{e_i}e_i=0,\quad\nabla_{e_1}e_2=a_2 e_1,\quad \nabla_{e_2}e_1=a_1 e_2\]
\[\Delta f=\theta\perp\extd f=-a_1\del_1 f-a_2\del_2 f,\quad \Delta e_1=2a_1e_2,\quad\Delta e_2=2a_2e_1,\quad  \sigma_{e_i}(e_i\tens e_i)={1\over 2}a_ie_i,\]
\[ \sigma_{e_1}(e_2\tens e_1)={1\over2}a_1e_2,\quad \sigma_{e_1}(e_2\tens e_2)=-a_2 e_1,\quad \sigma_{e_2}(e_1\tens e_1)=-a_1e_2,\quad \sigma_{e_2}(e_1\tens e_2)={1\over 2}a_2 e_1\]
and zero otherwise. These results come from  $\Vol\perp\theta=e_1 a_2-e_2 a_1$ so that
\[ \nabla_{e_i}e_j=-{1\over 2}(\eps_{ij}(e_1a_2-e_2a_1)-a_i e_j + e_i a_j),\quad  \sigma_{e_k}(e_i\tens e_j)={1\over 2}(\eps_{ki}e_{j'}\eps_{j'j}+e_k\delta_{ij})a_j\]
where $j'\ne j$.  Now since $\cj$ is invertible for $a,b\ne 0$ we look at the corresponding metric and connection:

Next, the parameter in $g$ up to overall normalisation is the one constant $a/b$, which also defines the function $\alpha$. We use $g$ to convert the form-connection coming from the cocycle to a connection $\Omega^1\to \Omega^1\tens_A\Omega^1$, which is straightforward noting that $e_i\alpha=-\alpha e_i$. This connection is torsion free since $\extd e_i=0$ and clearly $\wedge \nabla e_i=0$. One can check that it is also cotorsion-free. Here $\extd ({1\over a_1})=-{2\over a_1}e_2$ and $\extd ({1\over a_2})=-{2\over a_2}e_1$ from which
\[ coT_\nabla=\extd({2\over a_i} e_i)\tens e_i- {2\over a_i}e_i\wedge\nabla e_i=0\]
The connection is not, however, metric compatible as a bimodule connection. Its curvature is computed using $\extd \alpha^{\pm 1}=-2\alpha^{\pm 1}\theta$ and $[\alpha,\Vol]=0$. \endproof

Two byproducts of the cocycle construction were that we also have an `interior product' and a connection  on 2-forms, in the above example 
\[ j_{\Vol}(e_i)=-a_i\eps_{ij}e_j,\quad \nabla_{e_i}(\Vol)=0\]
from the formulae found for $\perp$. We also have a Hodge Laplacian on all degrees. The connection is a braided-derivation in that (\ref{wedgecompatg}) holds.

\subsection{Classical limit of the quantum Koszul formula}\label{Secpi}

As a small corollary of the quantum Koszul formula we apply it in the classical case of $A=C^\infty(M)$ for $M$ a Riemannian manifold with its classical exterior algebra $\Omega(M)$. However, we let $\perp= (\ ,\ ) + \pi$ on 1-forms instead of the obvious choice  $\perp_M=(\ ,\ )$, where $(\ ,\ )$ is the inverse metric and $\pi$ is an antisymmetric bivector field. 

First it can be shown that we can extend $\perp$ to higher forms by the same formula as in \cite{Ma:rec} (as recalled in the introduction) as an extended `inner product' but for the not-necessarily symmetric $(\ ,\ )+\pi$ on 1-forms. In particular, we have
\[ \omega\perp\eta=\omega\perp_M\eta+ (-1)^{|\omega|} L_{i_\pi}(\omega,\eta)\]
where $\perp_M$ is the usual extension of $(\ ,\ )$ and if $\pi=\pi_1\pi_2$ (sum of such terms understood) we define $i_\pi=i_{\pi_1}i_{\pi_2}$ as in \cite{Ma:rec} where $i$ along a vector field is the usual interior product. Thus 
\[ i_\pi(\omega\eta)=i_{\pi_1}(\pi_2(\omega)\eta-\omega i_{\pi_2}(\eta))=i_{\pi_1}(\eta)i_{\pi_2}(\omega)-i_{\pi_2}(\eta)i_{\pi_1}(\omega)\]
 if $\omega$ is a 1-form. Similarly, if $\delta_M$ is the usual Riemannian codifferential, we define
\[ \delta=\delta_M+ [ \extd, i_\pi] \]
and check
\begin{align*} \delta(a\omega)&=\delta_M(a\omega)+(\extd a)i_\pi(\omega)+a\extd i_\pi\omega -i_\pi((\extd a)\omega)-ai_\pi\extd \omega)\\
&=a\delta_M(\omega)+j_{\extd a}(\omega)+ a [\extd ,i_\pi]-L_{i_\pi}(\extd a,\omega)=a\delta(\omega)+\extd a\perp\omega\end{align*}
This is a special case (the classical limit) of \cite[Lem. 3.13]{Ma:rec}, which says that $\Delta, [[\ ,\ ]]$ are unchanged by adding the $\pi$ terms i.e. we still get the Riemannian Hodge Laplacian and Levi-Civita connection from our approach to the Koszul formula. 

In the extreme case we set $(\ ,\ )=0$ and $\delta_M=0$ so that $\perp=\pi$ on 1-forms. In this case our `connection' given by the cocycle obeys $\nabla_\omega(a\eta)=a\nabla_\omega\eta$ so $\nabla$ in this limit is a tensor.  

\section{Bicrossproduct model with $\alpha$-calculus}\label{Secalpha}

We let $A$ be the 2D bicrossproduct model spacetime algebra $A$ with generators $r,t$ and relations $[r,t]=\lambda r$ where $\lambda$ is an imaginary parameter. We consider the '$\alpha$-calculus' \cite{MaTao} given by commutation relations $[t, \extd r]=-\lambda\extd r$, $[t,\extd t]=\lambda\alpha\extd t$ and note that  in this case 
\begin{align*}[r^{\alpha},t]=&\lambda\alpha r^{\alpha}\\  
[t,\extd r^{\alpha}]=&\alpha[t,r^{\alpha-1}\extd r]=\alpha[t,r^{\alpha-1}]\extd r+ \alpha r^{\alpha-1}[t,\extd r]=-\lambda\alpha(\alpha-1)r^{\alpha-1}\extd r-\alpha r^{\alpha-1}\lambda\alpha\extd r\\
=&-\lambda\alpha^2r^{\alpha-1}\extd r=-\lambda\alpha\extd r^{\alpha}
\end{align*}
Thus if we set $r^{\alpha}\rightarrow r'$ and $\lambda\alpha\rightarrow\lambda'$ and then drop the prime notation, this is equivalent to setting $\alpha=1$ in our original differential algebra. Thus, as remarked in the introduction, we need only to consider this case.  We choose a central basis $e_1=r^{-1}\extd r$ and $e_2=r\extd t$. The exterior algebra is defined by $e_i^2=0$, $e_1e_2+e_2e_1=0$ and top form $e_1e_2=\extd r\extd t=\Vol$. We will see in this section how the quantum Koszul formula can be used to find the quantum Levi-Civita connection for any central quantum metric. We start by solving for $\perp$. 

\begin{lemma} Any matrix of constant entries $e_i\perp e_j=b_{ij}$ defines a solution of the 4-term relations with 
\[ e_i\perp\Vol=b_{ij}\eps_{jk}e_k,\quad \Vol\perp e_j=b_{ij}\eps_{ik}e_k,\quad \Vol\perp\Vol=\eps_{ij}b_{ij}\Vol\]
(sum of repeated indices). \end{lemma}
\proof  Because the $e_i$ are central we must have $a(e_i\perp e_j)=(ae_i\perp e_j=(e_i a)\perp e_j=e_i\perp (ae_j)=e_i\perp (e_ja)=(e_i\perp e_j)a$ for all $a\in A$, i.e. the $e_i\perp e_j$ must be in the centre of the algebra. In the polynomial setting the centre is the constants. The content of the 4-term relations in this case are otherwise exactly the same as the classical case and so it is not surprising that we find the same form as classically. We look at the 4-term relations for the various cases of 1-forms. If they all coincide, for example, 
\begin{align*}
-e_1e_1\perp e_1+(e_1\perp e_1)e_1=&e_1\perp e_1e_1+e_1(e_1\perp e_1)\ \Rightarrow\  b_{11}e_1=e_1b_{11}
\end{align*}
holds automatically as $e_1$ is central. Similarly for $e_2$. Next we have 
\begin{align*}
-e_1e_1\perp e_2+(e_1\perp e_1)e_2=&e_1\perp e_1e_2+e_1(e_1\perp e_2)\ \Rightarrow\  e_1\perp\Vol=&b_{11}e_2-e_1b_{12}
\end{align*}
\begin{align*}
-e_2e_1\perp e_1+(e_2\perp e_1)e_1=&e_2\perp e_1e_1+e_2(e_1\perp e_1)\ \Rightarrow\ \Vol\perp e_1=&e_2b_{11}-b_{21}e_1
\end{align*}
\begin{align*}
-e_1e_2\perp e_1+(e_1\perp e_2)e_1=&e_1\perp e_2e_1+e_1(e_2\perp e_1)\ \Rightarrow\ -\Vol\perp e_1+b_{12}e_1=&-e_1\perp\Vol+e_1b_{21}\\
\end{align*}
of which the first two are as stated and the last is then automatic. Similarly for $\Vol\perp e_2$ and $e_2\perp\Vol$ with the roles of $1,2$ interchanged. Finally, we look at the 4-term relations with $\omega=e_1,\eta=e_2,\zeta=\Vol$ which gives $\Vol\perp\Vol$ as stated. Other cases and other positions of $\Vol$ give nothing new. For example with $\eta=\Vol$ the 4-term relation requires
\[ (e_i\perp\Vol)e_j=-e_i(\Vol\perp e_j)\]
which holds for the solution found, again as is the case classically for $\perp$. \endproof

We also need to choose $\delta$ which we leave open and characterise by four functional parameters 
\begin{equation}\label{deltaab} \delta e_i=a_i \in A,\quad \delta\Vol= \sum_i b_i e_i,\quad b_i\in A\end{equation}
We similarly define matrices  by
\begin{equation}\label{jab} \cj_{e_i}(e_j)=(e_i,e_j)=g^{ij},\quad \cj_{\Vol}(e_i)=v^{ij}e_j\,\quad g^{ij},v^{ij}\in A\end{equation}
for the quantum metric/interior product that we  construct from $(\delta,\perp)$. 

\begin{proposition}\label{propalphaaibi}
\begin{enumerate}
\item For fixed $b_{ij}$, regular $\delta$ correspond to $a_i$ being at most linear in $t,r^{-1}$ and $b_i$ at most linear in $t,r$.\\
\item For all $g^{ij}, v^{ij}$ there exists a unique choice of $a_i,b_i$ up to constants $k_i, l_i$.\\
\item Non-singular $a_i,b_i$ correspond to $g=\frac{1}{2}(b+b^{T})$ to order $\lambda$, the symmetrisation of the matrix $b$, and $v^{i1}=-g^{i2}, \quad v^{i2}=g^{i1}$ to order $\lambda$ i.e. deforming the classical  interior product as a derivation. These hold exactly, not only to order $\lambda$, if and only if the $a_i,b_i$ are constants. \\
\item In the generic case where $|b|\ne 0$, $\delta^2$ is a left module map if and only 
\begin{align*} &(i)\quad\quad  g^{12}={1\over 2}(b_{12}+b_{21}),\quad g^{22}=b_{22},\quad v^{i1}=-g^{2i},\quad v^{i2}=g^{1i}\\
&(ii)\quad\quad l_1+k_2=b_{12},\quad l_2-k_1+b_{11}=0.\end{align*}
\item In the generic case where also $b_{11}b_{22}\ne b_{12}^2$, $\delta^2$ is a bimodule map if and only if in addition to (4),  $g^{ij}={1\over 2}(b_{ij}+b_{ji})$, or equivalently if and only if the $a_i,b_i$ are constants (related as in (4)(ii)).
\end{enumerate}
\end{proposition}
\proof (1) To apply Theorem 2.1 we need $\delta$ to be regular in the sense of a suitable bimodule map $\perp_R:\Omega\tens_A\Omega^1\to \Omega$. Since $e_i$ are central, if $\perp_R$ exists it must be given by
\[ e_i\perp_{R}\extd a=\extd a\perp e_i+[a,a_i]\]
and we take this as a definition extended as a bimodule map. It is well-defined since
\begin{align*} e_i\perp^R(a\extd b)&=e_i\perp^R(\extd(ab))-e_i\perp^R((\extd a)b)=\extd(ab)\perp e_i +[ab,a_i]-(\extd a\perp e_i + [a,a_i])b\\ 
&=\extd(ab)\perp e_i - \extd a\perp e_i b+ a[b,a_i]=\extd(ab)\perp e_i - ((\extd a)b)\perp e_i+ a[b,a_i]\\
&=a(\extd b\perp e_i+[b,a_i])=a(e_i\perp^R\extd b)=(e_ia)\perp^R\extd b. \end{align*}
We then compute $\cj_{e_i}(e_j)=\frac{1}{2}(e_i\perp e_j+e_i\perp_{R}e_j)$ which gives 
\[ (g^{ij})={1\over 2}\begin{pmatrix} 2b_{11}+r^{-1}[r,a_1] & b_{12}+b_{21}+r[t,a_1]\\ b_{12}+b_{21}+r^{-1}[r,a_2]\quad & 2 b_{22}+r[t,a_2]\end{pmatrix} \]
We also have 
\[ \Vol\perp_R e_1=(r^{-1}[r,b_1]-b_{12})e_1+(r^{-1}[r,b_1]+b_{11})e_2\]
\[\Vol\perp_R e_2=(r[t,b_1]-b_{22})e_1+(r[t,b_2]+b_{21})e_2\]
giving $\cj_{\Vol}(e_i)={1\over 2}(\Vol\perp e_i+\Vol\perp^R e_i)$ and therefore \[ \cj_{\Vol}(e_1)=-{1\over 2}(b_{12}+b_{21})e_1+b_{11}e_2+{1\over 2}r^{-1}[r,b_i]e_i\]
\[ \cj_{\Vol}(e_2)=-b_{22}e_1+{1\over 2}(b_{12}+b_{21})e_2+{1\over 2}r^{-1}[t,b_i]e_i.\]
We then want to invert these expressions to find the form of $a_i$ and $b_i$, ensuring that $g^{ij}$ and $v^{ij}$ remain constant parameters. We consider each component of the quantum metric separately. From the expression of $g^{11}$ we have that $a_1$ must be of the form $a_1=\frac{2}{\lambda}(g^{11}-b^{11})t+f(r)$ for some function $f$. Obtaining a particular $g^{12}$ then tells us that
\[ 2g^{12}-(b_{12}+b_{21})=r[t,a_1]=r[t,f(r)]=-\lambda r^2f'(r)\]
This has solution \[f(r)=\frac{2}{\lambda r}\left(g^{12}-\frac{1}{2}(b_{12}+b_{21})\right)+k_1\]
for some constant of integration $k_1$. This gives 
\begin{equation}\label{3a1}
a_1=\frac{2}{\lambda}\left((g^{11}-b_{11})t+\left(g^{12}-\frac{1}{2}(b_{12}+b_{21})\right)\frac{1}{r}\right)+k_1
\end{equation} Similarly, for $g^{21}$ we need $a_2=\frac{2}{\lambda}\left(g^{21}-\frac{1}{2}(b_{12}+b_{21})\right)t+g(r)$ for some function $g$. Then to obtain $g^{22}$ we need 
\[2g^{22}-2b_{22}=r[t,a_2]=r[t,g(r)]=-\lambda r^2g'(r)\] which has solution \[g(r)=\frac{2}{\lambda r}(g^{22}-b_{22})+k_2\] giving 
\begin{equation}\label{3a2}
a_2=\frac{2}{\lambda}\left(\left(g^{21}-\frac{1}{2}(b_{12}+b_{21})\right)t+(g^{22}-b_{22})\frac{1}{r}\right)+k_2
\end{equation}
We can see that $a_i$ has to be at most linear in $t$ and $r^{-1}$ in order for $g^{ij}$ to be constant and hence $\cj$ a bimodule map. For $b_i$ we consider
\[ (v^{ij})={1\over 2}\begin{pmatrix} -(b_{12}+b_{21)}+r^{-1}[r,b_1] & 2b_{11}+r^{-1}[r,b_2]\\ -2b_{22}+r^{-1}[t,b_1] & b_{12}+b_{21}+r^{-1}[t,b_2] \end{pmatrix} \] and we repeat the exact same process used to invert the $g^{ij}$. This gives 
\begin{equation}\label{3b1}
b_1=\frac{2}{\lambda}\left(\left(v^{11}+\frac{1}{2}(b_{12}+b_{21})\right)t-(v^{21}+b_{22})r\right)+l_1
\end{equation}
\begin{equation}\label{3b2}
b_2=\frac{2}{\lambda}\left((v^{12}-b_{11})t-\left(v^{22}-\frac{1}{2}(b_{12}+b_{21})\right)r\right)+l_2
\end{equation}
for constants of integration $l_i$. We can se that these are at most linear in $t,r$.

(2) The inverse metric coefficients $g^{ij}$ together with the coefficients $v^{ij}$ together form an 8-parameter space. Using a change of notation we can write 
\[a_i=\bar{a_i}t+\hat{a_i}r^{-1}+k_i, \quad b_i=\bar{b_i}t+\hat{b_i}r+l_i\]
which gives 12 parameters $\bar{a_i},\hat{a_i},\bar{b_i},\hat{b_i},k_i,l_i$. However, as the $a_i,b_i$ only ever appear as a commutation with either of the functions $r$ or $t$, the constants of integration do not effect the resulting values of $g^{ij}, v^{ij}$. Thus we are left with 8 genuine parameters, giving us a unique choice up to constants.

(3) Using the above notation, for the parameters to be non singular we need $\bar{a_i},\hat{a_i},\bar{b_i},\hat{b_i}$ to vanish to order $\lambda$. This happens precisely when we have the conditions stated. We assume that the constants $k_i,l_i$ are nonsingular as functions of $\lambda$, i.e. have a classical limit. 

(4) We compute
\begin{align*}\delta^2(f\Vol)&=\delta(f\delta\Vol+ \extd f\perp\Vol)=f\delta^2\Vol+\extd f\perp (b_i e_i)+\delta(\del^j f b_{jm}\eps_{mk}e_k)\\
&=f\delta^2\Vol+(\del^j f)(b_i b_{ji}+ b_{jm}\eps_{mk}a_k)+(\del^l\del^jf)b_{jm}\eps_{mk}b_{lk}\end{align*}
Requiring all but the first term to vanish for all $f$ gives
\[(\partial^1f)(b_{11}c_1+b_{12}c_2)+(\partial^2f)(b_{21}c_1+b_{22}c_2)+(\partial^2f\partial^1f)|b|-(\partial^1\partial^2f)|b|=0\]
for $c_1=b_1+a_2, c_2=b_2-a_1$. Here $\extd f=(\partial^1f)e_1+(\partial^2f)e_2$ define the partial derivatives. Since $r$ and $t$ generate the algebra, it suffices to require the above for $f=r,t$. These choices give 
\[b_{11}c_1+b_{12}c_2=0,\quad  b_{21}c_1+b_{22}c_2+|b|=0\]
with solution when $|b|\ne 0$,
\begin{equation}\label{alphadeltaleft} b_1+a_2=b_{12},\quad b_2-a_1=-b_{11}.\end{equation}
Inserting (\ref{3a1})-(\ref{3b2}) gives these in terms of the constant parameters as stated on looking at different powers of $t,r$. In principle there could be some further possibilities when $|b|=0$. 

(5) Since $\Vol$ is central, the condition for a $\delta^2$ to also be a right module map is that 
\begin{equation}\label{alphadeltabi} \delta^2\Vol=b_ia_i+(\del^j b_i)b_{ij}\end{equation}
is central (summations understood). To evaluate this we compute  $\extd b_1$ we find that $\del^1b_1=\frac{2}{\lambda r}(v^{21}+b_{22})=0$ (where we used the left-module map condition) and $\del^2b_1=\frac{2}{\lambda r}(g^{12}-g^{21})$. Similarly by considering $\extd b_2$ we find that $\del^1b_2=0$ and $\del^2b_2=\frac{2}{\lambda r}(g^{11}-b_{11})$. Substituting these expressions  $\delta^2(\Vol)$  we need
\[\frac{2t}{\lambda}((g^{12}-g^{21})(k_1-l_2)+(g^{11}-b_{11})(l_1+k_2))+\frac{2}{\lambda r}((g^{12}-g^{21})b_{12}+(g^{11}-b_{11})b_{22})+l_1k_1+l_2k_2\]
to be central. Applying the left-module conditions this becomes 
\[ (g^{12}-g^{21})b_{11}+(g^{11}-b_{11})b_{12}=0, \quad (g^{12}-g^{21})b_{12}+(g^{11}-b_{11})b_{22}=0\]
as our additional conditions to those of part (3). If $\det\begin{pmatrix} b_{11} & b_{12} \cr b_{12} & b_{12}\end{pmatrix}\neq0$ then this is equivalent to $g^{12}=g^{21}, g^{11}=b_{11}$ which given the results of part (3) is equivalent to $g=\frac{1}{2}(b+b^T)$ as matrices. There are some further exceptional cases where $\delta^2$ is a bimodule map and the above determinant vanishes.  Finally, we observe that the conditions displayed in (4)(i) and (5) of the proposition are equivalent to the conditions in part (3) for the $a_i,b_i$ to be constant. So apart from the exceptional cases, if $\delta^2$ is a left module map then it is a bimodule map if and only if the $a_i,b_i$ are constants $a_i=k_i, b_i=l_i$ (with $b_i$ determined from the $a_i$ by by (4)(ii)). From the above, its value is
\[\delta^2\Vol=l_1k_1+l_2k_2=b_{12}k_1-b_{11}k_2,\]
which includes zero as we can choose the remaining parameters freely.  \endproof

We are interested in obtaining $g^{ij}$  invertible with inverse $g_{ij}$ and metric $g=g_{ij}e_i\tens e_j$ central. This forces $g_{ij}$ to be constants (since the $e_i$ are central) and the coordinate algebra has a small centre. We may also want $g$ to be quantum symmetric in the sense $\wedge g=0$ which in our case just means $g_{ij}$ symmetric and hence in a real setting quantises AdS or dS geometry in 2D. Proposition~\ref{propalphaaibi} (3) says that this important case arises just from the assumption that $\delta$ has a classical limit. We also quantise the interior product $\cj$ in this case. We see that the same conclusion holds in (5) from requiring the algebraic property that $\delta^2$ is `strongly tensorial' in the sense of a bimodule map as in the classical case in \cite{Ma:rec}.

To complete the quantum geometry we proceed in the case $g$ invertible to construct the quantum connection associated to $(\delta,\perp)$ by the quantum Koszul formula in Theorem~2.1. We adopt the notations
\[ T_\nabla(e_i)=T_i\Vol,\quad coT_\nabla=\sum_i C_i\Vol\tens e_i;\quad T_i,C_i\in A\]
 to describe the resulting torsion and cotorsion. We will display the connection in the case where it has a classical limit, but the full expression can be found in the proof. We always take $\delta$ defined as they must be by $a_i, b_i$ in (\ref{3a1})-(\ref{3b2}) for given central invertible $g^{ij}$ and given $v^{ij}$. 
 
\begin{theorem}\label{thmalphalc} (1)  The resulting connection $\nabla$ depends only on the combinations $b_1+a_2,b_2-a_1$ and is non-singular if and only if the $\delta^2$ left module conditions (4)(i) in Proposition~\ref{propalphaaibi} hold to order $\lambda$. In this case the cotorsion and torsion and connection are
\[ C_1={1\over |g|}(b_1+a_2-b_{12}),\quad C_2={1\over |g|}(b_2-a_1+g^{11})\]
\[ T_1=\frac{1}{2}(g^{11}C_1+g^{12}C_2), \quad T_2=\frac{1}{2}(g^{12}C_1+g^{22}C_2)\]
\begin{align*}
\nabla e_1=&\frac{1}{2|g|}g^{12}(|g|C_1+2g^{12})e_1\tens e_1-\frac{1}{2|g|}g^{12}(2g^{11}-|g|C_2)e_1\tens e_2\\
&-\frac{1}{2|g|}g^{11}(|g|C_1+2g^{12})e_2\tens e_1+\frac{1}{2|g|}g^{11}(2g^{11}-|g|C_2)e_2\tens e_2\\
\nabla e_2=&\frac{1}{2|g|}g^{22}(|g|C_1+2g^{12})e_1\tens e_1+\frac{1}{2|g|}(g^{22}|g|C_2-2g^{12}g^{12})e_1\tens e_2\\
&-\frac{1}{2|g|}(g^{12}|g|C_1+2g^{11}g^{22})e_2\tens e_1-\frac{1}{2|g|}g^{12}(|g|C_2-2g^{11})e_2\tens e_2
\end{align*}
to order $\lambda$. 

(2) These formulae hold exactly, not only to order $\lambda$, if and only if the $\delta^2$ left module conditions (4)(i) in Proposition~\ref{propalphaaibi} hold. 

(3) The quantum connection in (2) is torsion free and metric compatible (i.e. $\nabla$ is a quantum Levi-Civita connection) if and only if the remaining $\delta^2$  bimodule map conditions displayed in (4)(ii) and (5) in Proposition~\ref{propalphaaibi} also hold (i.e. the $a_i,b_i$ are constants with $b_1+a_2=b_{12}$ and $b_2-a_1+b_{11}=0$). The associated braiding $\sigma(e_i\tens e_j)=e_j\tens e_i$. 
\end{theorem}
\proof (i) We compute the cocycle and hence 1-form covariant derivative from Theorem~2.1 as  \begin{align*}
\llbracket e_1,e_1\rrbracket=&0,\quad \llbracket e_1,e_2\rrbracket=(a_2+b_1-b_{12})e_1+(b_2-a_1+b_{11})e_2\\
\llbracket e_2,e_1\rrbracket=&-(a_2+b_1+b_{21})e_1+(a_1-b_2+b_{11})e_2,\quad \llbracket e_2,e_2\rrbracket=-2b_{22}e_1+(b_{12}+b_{21})e_2. 
\end{align*}
 We can also compute the braiding map $\sigma_\omega$ using the formula $\sigma(\eta\tens\zeta)=\cj_{\omega\eta}\zeta+\omega\cj_{\eta}\zeta$ and making use of $g^{ij}=\frac{1}{2}(b_{ij}+b_{ji})$  as, \begin{align*}
\sigma_{e_1}(e_1\tens e_1)=&\frac{1}{2}(2b_{11}+r^{-1}[r,a_1])e_1,\quad 
\sigma_{e_1}(e_1\tens e_2)=\frac{1}{2}(b_{12}+b_{21}+r[t,a_1])e_1\\
\sigma_{e_2}(e_1\tens e_1)=&\frac{1}{2}(b_{12}+b_{21}-r^{-1}[r,b_1])e_1+\frac{1}{2}r^{-1}([r,a_1]-[r,b_2])e_2\\
\sigma_{e_2}(e_1\tens e_2)=&\frac{1}{2}(2b_{22}-r^{-1}[t,b_1])e_1+\frac{1}{2}r^{-1}([t,a_1]-[t,b_2])e_2\\
\sigma_{e_1}(e_2\tens e_1)=&\frac{1}{2}r^{-1}([r,a_2]+[r,b_1])e_1+\frac{1}{2}(2b_{11}+r^{-1}[r,b_2])e_2\\
\sigma_{e_1}(e_2\tens e_2)=&\frac{1}{2}(r[t,a_2]+r^{-1}[t,b_1])e_1+\frac{1}{2}(b_{12}+b_{21}+r^{-1}[t,b_2])e_2\\
\sigma_{e_2}(e_2\tens e_1)=&\frac{1}{2}(b_{12}+b_{21}+r^{-1}[r,a_2])e_2,\quad \sigma_{e_2}(e_2\tens e_2)=\frac{1}{2}(2b_{22}+r[t,a_2])e_2
\end{align*}

(ii) We next define our abstract connection via the metric as $\nabla e_i=g^1\tens\nabla_{g^2}e_i$, where $g=g^1\tens g^2=g_{ij}e_i\tens e_j$ in terms of the inverse matrix $(g_{ij})$ which we write in terms of $(g^{ij})$ as usual. This gives
\begin{align*}
\nabla e_1=&\frac{1}{2}g^1\tens\llbracket g^2,e_1\rrbracket\\
=&\frac{1}{2|g|}g^{12}(b_1+a_2+b_{21})e_1\tens e_1-\frac{1}{2|g|}g^{12}(a_1-b_2+b_{11})e_1\tens e_2\\
&-\frac{1}{2|g|}g^{11}(b_1+a_2+b_{21})e_2\tens e_1+\frac{1}{2|g|}g^{11}(a_1-b_2+b_{11})e_2\tens e_2\\
\nabla e_2=&\frac{1}{2} g^1\tens\llbracket g^2,e_2\rrbracket\\
=&\frac{1}{2|g|}(g^{22}(b_1+a_2-b_{12})+2g^{12}b_{22})e_1\tens e_1+\frac{1}{2|g|}(g^{22}(b_2-a_1+b_{11})-g^{12}(b_{12}+b_{21}))e_1\tens e_2\\
&-\frac{1}{2|g|}(g^{21}(b_1+a_2-b_{12})+2g^{11}b_{22})e_2\tens e_1-\frac{1}{2|g|}(g^{21}(b_2-a_1+b_{11})-g^{11}(b_{12}+b_{21}))e_2\tens e_2
\end{align*}

(iii) For the torsion we compute  \begin{align*}
\wedge\nabla e_1-\extd e_1=&-\frac{1}{2|g|}g^{12}(a_1-b_2+b_{11})\Vol+\frac{1}{2|g|}g^{11}(b_1+a_2+b_{21})\Vol\\
=&\frac{1}{2|g|}(g^{11}(b_1+a_2)-g^{12}(a_1-b_2)+g^{11}b_{21}-g^{12}b_{11})\Vol\\
\wedge\nabla e_2-\extd e_2
=&\frac{1}{2|g|}(g^{22}(b_2-a_1+b_{11})-g^{12}(b_{12}+b_{21}))\Vol+\frac{1}{2|g|}(g^{21}(b_1+a_2-b_{12})+2g^{11}b_{22})\Vol-\Vol\\
=&\frac{1}{2|g|}(g^{22}(b_2-a_1)+g^{21}(b_1+a_2)+g^{22}b_{11}+2g^{11}b_{22}-g^{21}b_{12}-g^{12}(b_{12}+b_{21})-2|g|)\Vol
\end{align*} giving us
\begin{equation}\label{T1}T_1=\frac{1}{2|g|}(g^{11}(b_1+a_2)+g^{12}(b_2-a_1)+\frac{1}{2}g^{11}(b_{21}-b_{12}))
\end{equation}
\begin{equation}\label{T2}
T_2=\frac{1}{2|g|}(g^{22}(b_2-a_1)+g^{21}(b_1+a_2)+g^{11}b_{22}-g^{21}b_{12}).
\end{equation}

For cotorsion we compute 
\begin{align*}
(\extd\tens\id-\id\wedge\nabla) g=&-\extd\left(\frac{1}{|g|}g^{21}e_2\right)\tens e_1+\extd\left(\frac{1}{|g|}g^{11}e_2\right)\tens e_2 \\
&-\frac{1}{|g|}g^{22}e_1\nabla e_1+\frac{1}{|g|}g^{12}e_1\nabla e_2+\frac{1}{|g|}g^{21}e_2\nabla e_1-\frac{1}{|g|}g^{11}e_2\nabla e_2\\
=&\frac{1}{2|g|^2}g^{22}g^{11}(b_1+a_2+b_{21})\Vol\tens e_1-\frac{1}{2|g|^2}g^{22}g^{11}(a_1-b_2+b_{11})\Vol\tens e_2\\
&-\frac{1}{2|g|^2}(g^{12}g^{21}(b_1+a_2-b_{12})+2g^{11}g^{12}b_{22})\Vol\tens e_1\\
&-\frac{1}{2|g|^2}(g^{12}g^{21}(b_2-a_1+b_{11})-g^{11}g^{12}(b_{12}+b_{21}))\Vol\tens e_2\\
&-\frac{1}{2|g|^2}g^{21}g^{12}(b_1+a_2+b_{21})\Vol\tens e_1+\frac{1}{2|g|^2}g^{21}g^{12}(a_1-b_2+b_{11})\Vol\tens e_2\\
&+\frac{1}{2|g|^2}(g^{11}g^{22}(b_1+a_2-b_{12})+2g^{11}g^{12}b_{22})\Vol\tens e_1\\ 
&+\frac{1}{2|g|^2}(g^{11}g^{22}(b_2-a_1+b_{11})-g^{11}g^{12}(b_{12}+b_{21}))\Vol\tens e_2\\
&-\frac{1}{|g|}g^{21}\Vol\tens e_1+\frac{1}{|g|}g^{11}\Vol\tens e_2\\
=&\frac{1}{2|g|^2}(|g|(b_1+a_2+b_{21})+|g|(b_1+a_2-b_{12})-2|g|g^{21})\Vol\tens e_1\\
&+\frac{1}{2|g|^2}(|g|(b_2-a_1+b_{11})-|g|(a_1-b_2+b_{11})+2|g|g^{11})\Vol\tens e_2\\
=&\frac{1}{|g|}\left(b_1+a_2-g^{12}+\frac{1}{2}(b_{21}-b_{12})\right)\Vol\tens e_1+\frac{1}{|g|}(b_2-a_1+g^{11})\Vol\tens e_2
\end{align*}
giving us 
\begin{equation}\label{Ci} C_1=\frac{1}{|g|}\left(b_1+a_2-g^{12}+\frac{1}{2}(b_{21}-b_{12})\right),\quad C_2=\frac{1}{|g|}(b_2-a_1+g^{11})\end{equation}
in terms of $b_1+a_2,b_2-a_1$. These expressions for $T_i$ and $C_i$ are invertibly related to $\{b_1+a_2,b_2-a_1\}$, in particular 
\begin{equation}\label{bi} b_1+a_2=|g|C_1+g^{12}-\frac{1}{2}(b_{21}-b_{12}), \quad b_2-a_1=|g|C_2-g^{11},\end{equation}
which we then use in (\ref{T1}) and (\ref{T2}) to find $T_i$ in terms of $C_i$ as 
\[ T_1={1\over 2}(g^{11}C_1+g^{12}C_2),\quad T_2={1\over 2}(g^{21}C_1+g^{22}C_2)+ {1\over 2|g|}(g^{11}(g^{22}-b_{22})+g^{21}(g^{12}-{1\over 2}(b_{12}+b_{21}))).\] 
We can also use (\ref{bi}) to write the connection above in terms of $C_i$ to give
\begin{align*}
\nabla e_1=&\frac{1}{2|g|}g^{12}\left(|g|C_1+g^{21}+\frac{1}{2}(b_{12}+b_{21})\right)e_1\tens e_1-\frac{1}{2|g|}g^{12}(g^{11}+b_{11}-|g|C_2)e_1\tens e_2\\
&-\frac{1}{2|g|}g^{11}\left(|g|C_1+g^{21}+\frac{1}{2}(b_{12}+b_{21})\right)e_2\tens e_1+\frac{1}{2|g|}g^{11}(g^{11}+b_{11}-|g|C_2)e_2\tens e_2\\
\nabla e_2=&\frac{1}{2|g|}\left(g^{22}\left(|g|C_1+g^{21}-\frac{1}{2}(b_{12}+b_{21})\right)+2g^{12}b_{22}\right)e_1\tens e_1\\
&+\frac{1}{2|g|}(g^{22}(|g|C_2-g^{11}+b_{11})-g^{12}(b_{12}+b_{21}))e_1\tens e_2\\
&-\frac{1}{2|g|}\left(g^{21}\left(|g|C_1+g^{21}-\frac{1}{2}(b_{12}+b_{21})\right)+2g^{11}b_{22}\right)e_2\tens e_1\\
&-\frac{1}{2|g|}(g^{21}(|g|C_2-g^{11}+b_{11})-g^{11}(b_{12}+b_{21}))e_2\tens e_2
\end{align*}
This simplifies as stated when the $a_i,b_i$ are constant. 

(iv) We can see from (\ref{Ci}) and (\ref{3a1})-(\ref{3b2}) that $C_i$ and hence $T_i$ and $\nabla$ as found above are nonsingular if and only if 
\[ g^{12}={1\over 2}(b_{12}+b_{21}),\quad g^{22}=b_{22},\quad v^{i1}=-g^{2i},\quad v^{i2}=g^{1i}\]
hold to order $\lambda$ and in this case the torsion and cotorsion are related as stated to order $\lambda$ and $C_i$ as stated to order $\lambda$. These are exactly part (i) of the conditions for $\delta^2$ to be a left module map in Proposition~\ref{propalphaaibi} (4) (i.e. without the restriction on the $k_i,l_i$).  

(v) Finally suppose the conditions displayed in Proposition~\ref{propalphaaibi} (4)(i) so we are in the case of (iv). Then $C_1=0$ is exactly one of the conditions (\ref{alphadeltaleft}) in the proof of Proposition~\ref{propalphaaibi} (4), while $C_2=0$ becomes the other half of this if and only if $g^{11}=b_{11}$ which is  condition displayed in (5) in Proposition~\ref{propalphaaibi}. These combined assumptions are equivalent to $a_i,b_i$ constant with values shown by part (3) of Proposition~\ref{propalphaaibi}. 

We also find from our formulae for $\sigma_\omega$ that $\sigma_{e_i}(e_k\tens e_j)=g^{ij}e_k$. It then follows , as $g^{ij}$ is inverse to $g_{ij}$,  that  $\sigma(e_i\tens e_j)= g^1\tens\sigma_{ g^2}(e_i\tens e_j)=e_j\tens e_i$ (this does not mean it is the flip map on general elements, as it extends as a bimodule map). We then compute 
\[ \nabla_{e_k} g=g_{ij}(\nabla_{e_k} e_i\tens e_j+  e_i\tens \nabla_{e_k}e_j)=0.\]
on using the values of $\nabla$ found in (i).  \endproof

This is in line with the main result in \cite{Ma:rec} that the Levi-Civita connection arises in the classical case for a flat central extension with $\delta$ of classical type (such as $\delta^2=0$), but now in the quantum case provided only that $\delta$ has a classical limit. We can also compute the quantum curvature of the quantum connection given for non-singular $\delta$ in Theorem~\ref{thmalphalc}. As in Theorem~\ref{thmalphalc}, we continue here under the $\delta^2$ left module map assumption in part (2) of the theorem. The formula for curvature was recalled in Section~2.1.  In terms of cotorsion this amounts in our case to
\begin{align*}
R_{\nabla}(e_1)=&-\frac{1}{4|g|}(|g|^2C_2C_1+4g^{11}g^{12})\Vol\tens e_1+\frac{1}{4|g|}(4g^{11}g^{11}-|g|^2C_2^2)\Vol\tens e_2\\
R_{\nabla}(e_2)=&\frac{1}{4|g|}(|g|C_1(2g^{12}+|g|C_1)-2g^{22}(2g^{11}-|g|C_2))\Vol\tens e_1\\
&\frac{1}{4|g|}(2g^{11}-|g|C_2)(2g^{12}-|g|C_1)\Vol\tens e_2
\end{align*}
which is of particular interest when $C_i=0$ so that we have the quantum Levi-Civita connection by the theorem. 

In our above analysis we have concentrated on the connection acting on 1-forms, but the cocycle construction also gives it on forms of all degree. 
Continuing in our $\delta^2$ left module map assumption,  similar calculation from $2\nabla_{e_i}\Vol= \llbracket e_i,\Vol \rrbracket$ gives
\begin{equation}\label{nablaVol}\nabla_{e_1}\Vol={1\over 2}|g|C_2\Vol,\quad \nabla_{e_2}\Vol=-{1\over 2}|g|C_1\Vol\end{equation}
 using $\Vol\perp\Vol$ from Lemma~3.1. We see at the quantum Levi-Civita connection where $C_i=0$ that $\nabla_{e_i}\Vol=0$. We can also compute  
 \[ (\nabla_{e_i}e_1)e_2+ e_1\nabla_{e_i}e_2=\frac{|g|}{2}(C_2-C_1- g^{11}+b_{11})\Vol\]
which vanishes in the quantum Levi-Civita case. So these coincide, i.e. the derivation rule  (\ref{wedgecompatg}) holds for quantum Levi-Civita connection.

Another by-product of our theory is a Hodge-Laplacian given by $\Delta=\delta\extd+\extd\delta$, which we compute in the general case on some generators as \begin{align*}
\Delta(r)=&\delta\extd r=\delta(re_1)=\extd r\perp e_1+r\delta e_1=r(b_{11}+a_1)\\
\Delta(t)=&\delta\extd t=\delta(r^{-1}e_2)=-r^{-2}\extd r\perp e_2+r^{-1}a_2=r^{-1}(a_2-b_{12})\\
\Delta(e_1)=&\delta\extd e_1+\extd\delta e_1=\extd a_1=0\\
\Delta(e_2)=&\delta\extd e_2+\extd\delta e_2=\delta(\Vol)+\extd a_2=b_1e_1+b_2e_2\\
\Delta(\Vol)=&\delta\extd\Vol+\extd\delta\Vol=\extd(b_1e_1)+\extd(b_2e_2)=(a_1-g^{11})\Vol
\end{align*}

Finally, we might wonder if our choice of $\delta$ has a geometric picture in terms of the quantum Levi-Civita connection as is the case classically in the form of a divergence. We let $\ci_\eta(\omega)=\cj_\omega(\eta)$ be the left handed  `interior product' defined by $\cj$ and a candidate for the geometric codifferential that works at least in the classical case is $\ci\circ\nabla$. Recall that the connection depends only on the combinations $b_1+a_2,b_2-a_1$ so $\delta$ is not fixed for a particular choice of metric and connection. Proposition~\ref{propalphaaibi} part (2) tells us that this freedom corresponds to the value of $v^{ij}$ and we can fix it geometrically as follows.

\begin{lemma} For the quantum Levi-Civita connection we have $\delta=\ci\circ \nabla$ if and only if $a_1=g^{11}$, $a_2=g^{12}$ and $b_{12}=b_{21}$. In this case $b_i=0$ and $\delta^2=0$. 
\end{lemma} 
\proof We compute
\begin{align*}\ci\circ\nabla(e_1)=&\frac{1}{2|g|}g^{11}g^{12}(|g|C_1+2g^{12})-\frac{1}{2|g|}g^{12}g^{12}(2g^{11}-|g|C_2)\\
&-\frac{1}{2|g|}g^{11}g^{12}(|g|C_1+2g^{12})+\frac{1}{2|g|}g^{11}g^{22}(2g^{11}-|g|C_2)g^{22}\\
=&\frac{1}{2|g|}(|g|C_2(g^{12}g^{12}-g^{11}g^{22})+2g^{11}(g^{11}g^{22}-g^{12}g^{12}))=g^{11}-\frac{1}{2}|g|C_2\\
\ci\circ\nabla(e_2)=&\frac{1}{2|g|}g^{22}g^{11}(|g|C_1+2g^{12})+\frac{1}{2|g|}g^{12}(g^{22}|g|C_2-2g^{12}g^{12})\\
&-\frac{1}{2|g|}g^{12}(g^{12}|g|C_1+2g^{11}g^{22})-\frac{1}{2|g|}g^{12}g^{22}(|g|C_2-2g^{11})=g^{12}+\frac{1}{2}|g|C_1
\end{align*} 
so for quantum Levi-Civita connection $C_i=0$ we have the same as $\delta e_i$  if and only if $a_1=g^{11}$,  $a_2=g^{12}$.  This agrees with $\ci\circ\nabla\Vol$ precisely when $b_{12}=b_{21}$. \endproof 

This quantises the classical choice of $\delta$ within our 2-parameter moduli of values of $a_i$ that lead to the same quantum Levi-Civita connections. We also see that the geometric divergence $\delta$ requires $g^{ij}=b_{ij}$ or $(\ ,\ )=\perp$ which is the natural choice for the classical theory in \cite{Ma:rec}.

\section{Bicrossproduct model with its standard differential calculus}\label{Secbeta}

The same quantum spacetime $A$ as in the previous section has another family of calculi, the $\beta$ calculus, for which the standard case ($\beta=1$) is  given by commutation relations 
\[ [r,\extd t]=\lambda\extd r,  \quad [t,\extd t]=\lambda\extd t,\quad [r,\extd r]=0,\quad [t,\extd r]=0\]
The general $\beta$ case is significantly more complicated but not expected to be fundamentally different in view of related work such as \cite{MaTao}. This time a central basis is 
\[ e_1=\extd r,\quad e_2=\nu=r\extd t-t\extd r.\]
and the canonical exterior algebra here obeys $e_1^2=0$, $e_2^2=-\lambda\Vol$ and $e_1e_2+e_2e_1=0$, with top form $\Vol=e_1e_2$. 

\begin{lemma} \label{perpbic} The general solution to the 4-term relations on this exterior algebra when $\lambda\ne 0$ has the form 
\[ \Vol\perp e_i=-e_i\perp\Vol= b e_i,\quad e_1\perp e_1=0,\quad e_1\perp e_2=b ,\quad e_2\perp e_1=-b,\quad e_2\perp e_2=-\lambda b\]
\[ \Vol\perp\Vol=2b\Vol\]
for some constant parameter $b$. 
\end{lemma}
\proof To start with we set $e_i\perp e_j=b_{ij}$ and require that $\perp$ is a bimodule map, which as in Section~3 forces the $b_{ij}$ to be constants.  The 4-term relation on $e_1,e_1,e_1$ gives that $[b_{11},e_1]=0$  as in the proof of Lemma~3.1, which is automatic. On $e_1,e_1,e_2$ we have
\begin{align*}
-e_1e_1\perp e_2+(e_1\perp e_1)e_2=&e_1\perp e_1e_2+e_1(e_1\perp e_2)\ \Rightarrow\ b_{11}e_2=e_1\perp\Vol+e_1b_{12}
\end{align*}
Next, $e_2,e_1,e_1$ gives
\begin{align*}
-e_2e_1\perp e_1+(e_2\perp e_1)e_1=&e_2\perp e_1e_1+e_2(e_1\perp e_1)\ \Rightarrow\  \Vol\perp e_1+b_{21}e_1=e_2b_{11}
\end{align*}
The $e_1,e_2,e_1$ equation is automatic while  $e_1,e_2,e_2$ gives
\begin{align*}
-e_1e_2\perp e_2+(e_1\perp e_2)e_2=&e_1\perp e_2e_2+e_1(e_2\perp e_2) \Rightarrow\ -\Vol\perp e_2+b_{12}e_2=-\lambda e_1\perp\Vol+e_1b_{22}
\end{align*}
which in view of our previous values we write as 
\[ \Vol\perp e_2=(b_{12}+\lambda b_{11})e_2-(b_{22}+\lambda b_{12})e_1\]
Similarly on $e_2,e_1,e_2$ we have
\begin{align*}
-e_2e_1\perp e_2+(e_2\perp e_1)e_2=&e_2\perp e_1e_2+e_2(e_1\perp e_2)\ \Rightarrow\ \Vol\perp e_2+b_{21}e_2=e_2\perp\Vol+e_2b_{12}
\end{align*}
which we write as
\[ e_2\perp\Vol=(b_{21}+\lambda b_{11})e_2-(b_{22}+\lambda b_{12})e_1.\]
On $e_2,e_2,e_1$ we have
\begin{align*}
-e_2e_2\perp e_1+(e_2\perp e_2)e_1=&e_2\perp e_2e_1+e_2(e_2\perp e_1)\ \Rightarrow\ \lambda\Vol\perp e_1+b_{22}e_1=-e_2\perp\Vol+e_2b_{21}
\end{align*}
which we write as
\[ e_2\perp\Vol=(b_{21}-\lambda b_{11})e_2-(b_{22}-\lambda b_{21})e_1.\]
Comparing the two different values we have for $e_2\perp\Vol$ implies for $\lambda\ne 0$ that $b_{11}=0$ and $b_{12}=-b_{21}$. Finally, the 4-term relation on $e_2,e_2,e_2$ gives us 
\begin{align*}
-e_2e_2\perp e_2+(e_2\perp e_2)e_2=&e_2\perp e_2e_+e_2(e_2\perp e_)\ \Rightarrow\ \lambda\Vol\perp e_2+b_{22}e_2=-\lambda e_2\perp\Vol+e_2b_{22}
\end{align*}
which implies that $-\Vol\perp e_2=e_2\perp\Vol$ provided $\lambda\ne 0$. Comparing the values already obtained for these, we deduce that  $b_{22}=-\lambda b_{12}$. This gives the stated form with $b_{12}=b$. We also look at the 4-term relations with one of the forms being $\Vol$ to obtain the value shown. \endproof

This is already far from the classical case as the classical limit of $e_i\perp e_j$ is antisymmetric. We can still proceed to see what kinds of metrics and connections can be obtained by the quantum Koszul formula. As before, we take a general form of $\delta$ as in Section~3 namely $\delta e_i=a_i$ and $\delta\Vol=\sum_i b_ie_i$ for $a_i,b_i\in A$. 

\begin{proposition}\label{propbetaaibi}
\begin{enumerate}
\item For fixed constant parameter $b$, regular $\delta$ correspond to $a_i$ and $b_i$ each being at most linear in $\frac{t}{r}, \frac{1}{r}$.\\
\item For all $g^{ij},v^{ij}$ there exists a unique choice of $a_i,b_i$ up to constants $k_i,l_i$.\\
\item Nonsingular $a_i,b_i$ correspond to $g^{ij}=v^{ij}=0$ to order $\lambda$.\\
\item $\delta^2$ is a left module map if and only if 
\[ v^{i2}=g^{1i},\quad v^{i1}=-g^{2i}-\lambda g^{1i},\quad  l_1+k_2+\lambda k_1=0, \quad l_2-k_1=0. \]
\item $\delta^2$ is a bimdolue map if in addition 
\[ g^{11}=0, \quad g^{12}+g^{21}=0, \quad k_1=0\]
where the last two apply in the generic case of $b\ne 2g^{12}$.
\end{enumerate}
\end{proposition}

\proof The proof that $\delta$ is regular and that $\perp$ is a bimodule map is exactly the same as in the proof of Proposition~\ref{propalphaaibi}. Here we again use the equation for $\perp_R$, but for this calculus we find \begin{align*}
e_1\perp_R e_1=&[r,a_1],\quad e_1\perp_R e_2=[t,a_1]r-[r,a_1](t-\lambda)-b\\
e_2\perp_R e_1=&[r,a_2]+b,\quad e_2\perp_R e_2=[t,a_2]r-[r,a_2](t-\lambda)-\lambda b\\
\Vol\perp_R e_1=&([r,b_1]-b)e_1+[r,b_2]e_2\\
 \Vol\perp_R e_2=&([t,b_1]r-[r,b_1](t-\lambda))e_1+([t,b_2]r-[r,b_2](t-\lambda)-b)e_2
\end{align*} We then set $\cj_{e_i}(e_j)=g^{ij}$ for our quantum metric, $\cj_{\Vol}(e_i)=v^{ij}e_j$ and use our known data for $\perp$. This gives 
\[ (g^{ij})={1\over 2}\begin{pmatrix} [r,a_1] &  [t,a_1]r-[r,a_1](t-\lambda) \\ [r,a_2]\quad & [t,a_2]r-[r,a_2](t-\lambda)-2\lambda b\end{pmatrix} \]
\[ \cj_\Vol(e_1)=\frac{1}{2}([r,b_1]e_1+[r,b_2]e_2)\]
\[ \cj_\Vol(e_2)=\frac{1}{2}(([t,b_1]r-[r,b_1](t-\lambda))e_1+([t,b_2]r-[r,b_2](t-\lambda))e_2)\]
As before, we then want to invert this relationship and solve for $a_i$ and $b_i$ in such a way that $g^{ij}$ and $v^{ij}$ are constants (numerical parameters). We consider each component of the quantum metric separately. From the expression for $g^{11}$ we have that $a_1$ must be of the form $a_1=\frac{2tg^{11}}{\lambda r}+f(r)$ for some function $f$. We note here that the notation here $\frac{t}{r}$ always to be read $\frac{1}{r}\cdot t$ in our calculus. Obtaining a particular $g^{12}$ then tells us that 
\[  2g^{12}+2g^{11}(t-\lambda)= [t,a_1]r= [t,\frac{1}{r}]\frac{2tg^{11}}{\lambda}r+[t,f(r)]r\\
=\frac{2tg^{11}}{r}r+[t,f(r)]r\]
\[ =2g^{11}\frac{1}{r}(rt-\lambda r)+[t,f(r)]r
=2g^{11}(t-\lambda)+[t,f(r)]r\]
using the algebra commutation relations. Comparing the two sides, we see that $[t,f(r)]r=2g^{12}$ or $-\lambda f'(r)=\frac{2g^{12}}{r^2}$, which has soltution 
\[ f(r)=\frac{2g^{12}}{\lambda r}+k_1\]
for some constant of integration $k_1$.  This gives the form of $a_1$, namely 
\begin{equation}\label{4a1}
a_1=\frac{2}{\lambda r}(g^{11}t+g^{12})+k_1
\end{equation}
Similarly, for $g^{21}$ we need $a_2=\frac{2g^{21}t}{\lambda r}+g(r)$ for some function $g$. Then to obtain a particular $g^{22}$ we need
\[ 2g^{22}+2\lambda b+2g^{21}(t-\lambda) =  [t,a_2]r=[t,\frac{1}{r}]\frac{2g^{21}t}{\lambda}r+[t,g(r)]r=2g^{21}(t-\lambda)+[t,g(r)]r\]
Comparing the two sides we need deduce $[t,g(r)]r=2g^{22}+2\lambda b$ or $g'(r)=-\frac{2g^{22}}{\lambda r^2}-\frac{2b}{r^2}$ 
with solution
\[ g(r)=\frac{2g^{22}}{\lambda r}+\frac{2b}{r}+k_2\]
giving the form of $a_2$,
\begin{equation}\label{4a2}
a_2=\frac{2}{\lambda r}(g^{21}t+g^{22})+\frac{2b}{r}+k_2
\end{equation}
We can see that $a_i$ has to be at most linear in $\frac{t}{r},\frac{1}{r}$. For $b_i$ we consider 
\[ (v^{ij})={1\over 2}\begin{pmatrix} [r,b_1] &  [r,b_2] \\ [t,b_1]r-[r,b_1](t-\lambda) \quad & ([t,b_2]r-[r,b_2](t-\lambda), \end{pmatrix} \]
and we use the same process we used to invert for $g^{ij}$. This gives
\begin{equation}\label{4bi}
b_1=\frac{2}{\lambda r}(v^{11}t+v^{21})+l_1,\quad 
b_2=\frac{2}{\lambda r}(v^{12}t+v^{22})+l_2
\end{equation}
for constants of integration $l_i$. Again, we can observe that these are at most linear in $\frac{t}{r},\frac{1}{r}$. 

Parts (2) and (3) are obtained by solving as in Proposition~\ref{propalphaaibi}, only replacing $t$ with $\frac{t}{r}$ and for $b_i$ the constant $\hat{b_i}$ is used for the coefficient of $\frac{1}{r}$ as opposed to $r$. The form of (\ref{4a1})-(\ref{4bi}) tells us the conditions for $a_i,b_i$ to be non-singular in $\lambda$ assuming the $k_i,l_i$ are. For part (4) we compute 
\begin{align*}
\delta^2(f\Vol)=&\delta(f\delta\Vol+\extd f\perp\Vol)=f\delta^2\Vol+\extd f\perp b_ie_i-b \delta((\del^if)e_i)\\
=&f\delta^2\Vol+(\del^jf)(e_j\perp e_i)b_i-b(\del^if)a_i-(\del^l\del^if)b(e_l\perp e_i)
\end{align*}
Requiring all but the first term to vanish for all $f$ gives the condition
\[(\del^1f)b(b_2-a_1)-b(\del^2f)(b_1+a_2+\lambda b_2)+b^2(\del^2\del^1f)-b^2(\del^1\del^2f)+b^2\lambda(\del^2\del^2f)=0\]
where the partial derivatives are defined by $\extd$ in our basis $\{e_i\}$ as usual.  Again, since $r$ and $t$ generate the algebra it suffices to consider $f=t,r$ which respectively give the two conditions 
\begin{equation}\label{delta2leftab}b_1+a_2+\lambda b_2-\frac{2b}{r}=0,\quad b_2-a_1=0.\end{equation} We then use (\ref{4a1})-(\ref{4bi}) and consider different powers of $r,t$ to obtain the displayed equations in terms of $g^{ij},v^{ij},l_i,k_i$, using the first pair to present the 2nd pair as values of $l_i,v^{i1}$. 

For part (5), since $\Vol$ is central, the additional condition for a right and hence bi-module map is that 
\begin{equation}\label{delta2centralab} \delta^2\Vol=b_ia_i+b(  \del^1 -\lambda  \del^2)b_2 - b \del^2 b_1\end{equation}
is central. From $\extd b_i$ we find 
\[ \del^1b_1=-\frac{2v^{21}}{\lambda r^2}, \quad \del^2b_1=\frac{2v^{11}}{\lambda r^2}, \quad \del^1b_2=-\frac{2v^{22}}{\lambda r^2}, \quad \del^2b_2=\frac{2v^{12}}{\lambda r^2}\]
hence for (\ref{delta2centralab}) to be central we need 
\begin{align*}
&\frac{4}{\lambda^2r^2}(v^{11}t+v^{21})(g^{11}t+g^{12})+\frac{4}{\lambda r^2}v^{11}(g^{11}t+g^{12})+\frac{4}{\lambda^2r^2}(v^{12}t+v^{22})(g^{21}t+g^{22})\\
&+\frac{4}{\lambda r^2}v^{12}(g^{21}t+g^{22})+\frac{4b}{\lambda r^2}(v^{12}t+v^{22})+\frac{4b}{r^2}v^{12}-\frac{2bv^{22}}{\lambda r^2}-\frac{2bv^{12}}{r^2}-\frac{2bv^{11}}{\lambda r^2}\\
&+\frac{2k_1}{\lambda r}(v^{11}t+v^{21})+\frac{2l_1}{\lambda r}(g^{11}t+g^{12})+\frac{2k_2}{\lambda r}(v^{12}t+v^{22})+\frac{2l_2}{\lambda r}(g^{21}t+g^{22})+\frac{2bl_2}{r}\\
\end{align*}
to be central. At our level of polynomials in $t,r,r^{-1}$, we require the expression itself to vanish (leaving a constant $\delta^2\Vol=l_ik_i=-\lambda k^2$). Applying the left module map condition and collecting terms of order $\frac{1}{r^2}$ we have 
\begin{align*}
&-\frac{4}{\lambda}g^{11}g^{11}t^2-\frac{4}{\lambda}g^{12}g^{21}-\frac{4}{\lambda}g^{11}g^{12}t-\frac{4}{\lambda}g^{11}g^{12}t-\frac{4}{\lambda}g^{12}g^{12}\\
&-4g^{11}g^{11}t-4g^{11}g^{12}+4bg^{11}+\frac{4}{\lambda}g^{11}g^{22}+\frac{4b}{\lambda}g^{11}t+\frac{2b}{\lambda}g^{12}+\frac{2b}{\lambda}g^{21}
\end{align*}
For $t^2$ term to vanish we require $g^{11}=0$ as stated. Given this, the other terms vanish if and only $(g^{12}+g^{21})(b-2g^{12})=0$. We then examine the coefficient of the $\frac{1}{r}$ terms of our original expression which we again need to vanish in order to be central. Assuming we have $\delta^2$ a left module map, we are left with $2(b-2g^{12})k_1=0$.  For generic $b$ this means $g^{12}+g^{21}=0$ and $k_1^2=0$, and hence $\delta^2\Vol=0$. 
\endproof

Unlike Section~3, we see that we cannot usefully take $a_i, b_i$ and hence $\delta$ to be nonsingular in the sense of having a classical limit, if we want non-zero $g,\cj$ in the classical limit. However, we can still explore the resulting quantum geometry and ask for $\delta^2$ to be tensorial (at least a left module map). For a fixed $g^{ij},v^{ij}$ we the $a_i,b_i$ are uniquely defined according to the above by  (\ref{4a1})-(\ref{4bi}) up to free parameters $k_i,l_i$.  These play the role of the constant values of $a_i, b_i$ in Section~3 and do not affect the metric or $\cj$ but do affect the central extension cocycle and bimodule connection coming out of the quantum Koszul formula for our choice of $\delta$. We will study the quantum connection through its torsion and cotorsion coefficients $T_i, C_i$ defined as before. We let 
\[ |g|=\det(g^{ij}),\quad |g|_\lambda=|g|-{1\over 2}\lambda^2 (g^{11})^2.\]

\begin{lemma}\label{lemTC} The connection from the quantum Koszul formula for any fixed $g^{ij}$ and $v^{ij}$ has torsion and cotorsion
\begin{align*}
T_1=&\frac{1}{2|g|}\left(g^{11}\left(a_2+b_1-\frac{2b}{r}\right)-(g^{12}+\lambda g^{11})(a_1-b_2)\right)\\
T_2=&\frac{1}{2|g|}\left((b_2-a_1)(g^{22}+\lambda g^{21})+\lambda b_1g^{11}+\lambda b_2(g^{12}+\lambda g^{11})+g^{21}\left(a_2+b_1-\frac{2b}{r}\right)-\frac{4|g|}{r}\right)\\
C_1=&\frac{1}{|g|}\left(a_2+\frac{|g|_{\lambda}}{|g|}b_1-\frac{2b}{r}-\frac{2g^{21}}{r}\right)\\
C_2=&\frac{1}{|g|}\left(\frac{|g|_{\lambda}}{|g|}b_2-a_1+\frac{2g^{11}}{r}\right)
\end{align*}
where $a_i,b_i$ are given in terms of the parameters $k_i,l_i$ by (\ref{4a1})-(\ref{4bi}).
\end{lemma}
\proof (i) The covariant derivative along 1-forms is given by $\nabla_{e_i}e_j={1\over 2}[[e_i,e_j]]$ where the cocycle data in Theorem~\ref{thm4term} comes out as 
 \begin{align*}
\llbracket e_1,e_1\rrbracket=&0,\quad  \llbracket e_1,e_2\rrbracket=\left(a_2+b_1-\frac{2}{r}b\right)e_1+(b_2-a_1)e_2\\
\llbracket e_2,e_1\rrbracket=&\left(\frac{2}{r}b-a_2-b_1\right)e_1+(a_1-b_2)e_2,\quad  \llbracket e_2,e_2\rrbracket=-\lambda(b_1e_1+b_2e_2)
\end{align*} while the generalised braiding from $\sigma_{\omega}(\eta\tens\zeta)=\cj_{\omega\eta}\zeta+\omega\cj_{\eta}\zeta$ comes out as 
\begin{align*}
\sigma_{e_1}(e_1\tens e_1)=&\frac{1}{2}[r,a_1]e_1\\
\sigma_{e_1}(e_1\tens e_2)=&\frac{1}{2}([t,a_1]r-[r,a_1](t-\lambda))e_1\\
\sigma_{e_2}(e_1\tens e_1)=&\frac{1}{2}([r,a_1]-[r,b_1])e_1+\frac{1}{2}[r,b_2]e_2\\
\sigma_{e_2}(e_1\tens e_2)=&+\frac{1}{2}([r,b_1](t-\lambda)-[t,b_1]r)e_1+\frac{1}{2}([t,a_1]r-[r,a_1](t-\lambda)-[t,b_2]r+[r,b_2](t-\lambda))e_2\\
\sigma_{e_1}(e_2\tens e_1)=&\frac{1}{2}([r,b_1]+[r,a_2])e_1+\frac{1}{2}[r,b_2]e_2\\
\sigma_{e_1}(e_2\tens e_2)=&\frac{1}{2}([t,b_1]r-[r,b_1](t-\lambda)+[t,a_2]r-[r,a_2](t-\lambda)-2\lambda b)e_1+\frac{1}{2}([t,b_2]r-[r,b_2](t-\lambda))e_2\\
\sigma_{e_2}(e_2\tens e_1)=&-\frac{\lambda}{2}[r,b_1]e_1+\frac{1}{2}([r,a_2]-[r,b_2])e_2\\
\sigma_{e_2}(e_2\tens e_2)=&\frac{\lambda}{2}([r,b_1](t-\lambda)-[t,b_1]r)e_1+\frac{1}{2}([t,a_2]r-[r,a_2](t-\lambda)-2\lambda b-\lambda[t,b_2]r+\lambda[r,b_2](t-\lambda))e_2
\end{align*} 

(ii) The abstract connection is $\nabla e_i= g^1\tens\nabla_{ g^2}e_i$, where $ g= g^1\tens g^2=g_{ij}e_i\tens e_j$ is the metric with $(g_{ij})$ inverse to $(g^{ij})$. This comes out as
 \begin{align*}
\nabla e_1=&\frac{1}{2} g^1\tens\llbracket g^2,e_1\rrbracket\\
=&-\frac{1}{2|g|}g^{12}\left(\frac{2b}{r}-a_2-b_1\right)e_1\tens e_1-\frac{1}{2|g|}g^{12}(a_1-b_2)e_1\tens e_2\\
&+\frac{1}{2|g|}g^{11}\left(\frac{2b}{r}-a_2-b_1\right)e_2\tens e_1+\frac{1}{2|g|}g^{11}(a_1-b_2)e_2\tens e_2
\end{align*}
\begin{align*}
\nabla e_2=&\frac{1}{2} g^1\tens\llbracket g^2,e_2\rrbracket\\
=&\frac{1}{2|g|}\left(g^{22}\left(a_2+b_1-\frac{2b}{r}\right)+\lambda g^{12}b_1\right)e_1\tens e_1+\frac{1}{2|g|}(g^{22}(b_2-a_1)+\lambda g^{12}b_2)e_1\tens e_2\\
&-\frac{1}{2|g|}\left(g^{21}\left(a_2+b_1-\frac{2b}{r}\right)+\lambda g^{11}b_1\right)e_2\tens e_1-\frac{1}{2|g|}(g^{21}(b_2-a_1)+\lambda g^{11}b_2)e_2\tens e_2
\end{align*}

(iii)We can now compute the associated torsion $T_\nabla$ as \begin{align*}
\wedge\nabla e_1-\extd e_1
=&-\frac{1}{2|g|}g^{12}\left(a_1-b_2\right)\Vol+\frac{1}{2|g|}g^{11}\left(a_2+b_1-\frac{2b}{r}\right)\Vol-\frac{\lambda}{2|g|}g^{11}\left(a_1-b_2\right)\Vol\\
\wedge\nabla e_2-\extd e_2
=&\frac{1}{2|g|}\left(g^{22}\left(b_2-a_1\right)+\lambda g^{12}b_2\right)\Vol+\frac{1}{2|g|}\left(g^{21}\left(a_2+b_1-\frac{2b}{r}\right)+\lambda g^{11}b_1\right)\Vol\\
&+\frac{\lambda}{2|g|}\left(g^{21}\left(b_2-a_1\right)+\lambda g^{11}b_2\right)\Vol-\frac{2}{r}\Vol\\
\end{align*} 
from which we read off the values of $T_i$ as stated.  For the cotorsion, need \begin{align*}
coT_{\nabla}=&(\extd\tens\id-\id\wedge\nabla)\left(\frac{1}{|g|}(g^{22}e_1\tens e_1-g^{12}e_1\tens e_2-g^{21}e_2\tens e_1+g^{11}e_2\tens e_2)\right)
\end{align*} which we examine term by term:
\begin{align*}
(\extd\tens\id&-\id\wedge\nabla)\left(\frac{1}{|g|}g^{22}e_1\tens e_1\right)
=-\frac{1}{|g|}g^{22}e_1\nabla e_1\\
=&\frac{1}{2|g|^2}g^{22}g^{11}\left(a_2+b_1-\frac{2b}{r}\right)\Vol\tens e_1-\frac{1}{2|g|^2}g^{22}g^{11}\left(a_1-b_2\right)\Vol\tens e_2
\end{align*}
\begin{align*}
(\extd\tens\id&-\id\wedge\nabla)\left(-\frac{1}{|g|}g^{12}e_1\tens e_2\right)
=\frac{1}{|g|}g^{12}e_1\nabla e_2\\
=&-\frac{1}{2|g|^2}g^{12}g^{21}\left(a_2+b_1-\frac{2b}{r}\right)\Vol\tens e_1-\frac{\lambda}{2|g|^2}g^{12}g^{11}c\Vol\tens e_1\\
&-\frac{1}{2|g|^2}g^{12}g^{21}\left(b_2-a_1\right)\Vol\tens e_2-\frac{\lambda}{2|g|^2}g^{12}g^{11}d\Vol\tens e_2
\end{align*}
\begin{align*}
(\extd\tens\id&-\id\wedge\nabla)\left(-\frac{1}{|g|}g^{21}e_2\tens e_1\right)
=-\frac{1}{|g|}g^{21}\frac{2}{r}\Vol\tens e_1+\frac{1}{|g|}g^{21}e_2\nabla e_1\\
=&-\frac{1}{2|g|^2}g^{21}\frac{4|g|}{r}\Vol\tens e_1-\frac{1}{2|g|^2}g^{21}g^{12}\left(a_2+b_1-\frac{2b}{r}\right)\Vol\tens e_1\\
&+\frac{1}{2|g|^2}g^{21}g^{12}\left(a_1-b_2\right)\Vol\tens e_2+\frac{\lambda}{2|g|^2}g^{21}g^{11}\left(a_2+b_1-\frac{2b}{r}\right)\Vol\tens e_1\\
&-\frac{\lambda}{2|g|^2}g^{21}g^{11}\left(a_1-b_2\right)\Vol\tens e_2
\end{align*}
\begin{align*}
(\extd\tens\id&-\id\wedge\nabla)\left(\frac{1}{|g|}g^{11}e_2\tens e_2\right)
=\frac{1}{|g|}g^{11}\frac{2}{r}\Vol\tens e_2-\frac{1}{|g|}g^{11}e_2\nabla e_2\\
=&\frac{1}{2|g|^2}g^{11}\frac{4|g|}{r}\Vol\tens e_2+\frac{\lambda}{2|g|^2}g^{11}g^{12}c\Vol\tens e_1+\frac{1}{2|g|^2}g^{11}g^{22}\left(a_2+b_1-\frac{2b}{r}\right)\Vol\tens e_1\\
&+\frac{1}{2|g|^2}g^{11}g^{22}\left(b_2-a_1\right)\Vol\tens e_2+\frac{\lambda}{2|g|^2}g^{11}g^{12}d\Vol\tens e_2\\
&-\frac{\lambda}{2|g|^2}g^{11}g^{21}\left(a_2+b_1-\frac{2b}{r}\right)\Vol\tens e_1-\frac{\lambda^2}{2|g|^2}(g^{11})^2b_1\Vol\tens e_1\\
&-\frac{\lambda}{2|g|^2}g^{11}g^{21}\left(b_2-a_1\right)\Vol\tens e_2-\frac{\lambda^2}{2|g|^2}(g^{11})^2b_2\Vol\tens e_2
\end{align*} Collecting like coefficients of $\Vol\tens e_i$ and simplifying gives the coefficients $C_i$.  \endproof

We now want to look carefully at the classical limit and, knowing from Proposition~\ref{propbetaaibi} that $a_i,b_i$ will have to be singular for a nonzero geometry, we write them in terms of new parameters where we factor out an order $1/\lambda$ singularity, thus 
\begin{equation}\label{deltatilde}\delta(e_1)=\frac{\tilde{a_1}}{\lambda}, \quad \delta(e_2)=\frac{\tilde{a_2}}{\lambda}+\frac{2b}{r},\quad\delta\Vol={1\over\lambda}(\tilde b_1 e_1+\tilde be_2 e_2);\end{equation}
\begin{equation}\label{tildeai}
\tilde{a_1}=\frac{2}{r}(g^{11}t+g^{12})+\tilde{k_1},\quad \tilde{a_2}=\frac{2}{r}(g^{21}t+g^{22})+\tilde{k_2}
\end{equation}
\begin{equation}\label{tildebi}
\tilde{b_1}=\frac{2}{r}(v^{11}t+v^{21})+\tilde{l_1},\quad \tilde{b_2}=\frac{2}{r}(v^{12}t+v^{22})+\tilde{l_2}
\end{equation}
as the general form of regular $\delta$ in terms of rescaled constant parameters  $\tilde{k_i}=\lambda k_i,$ $\tilde{l_i}=\lambda l_i$. This is equivalent to our previous $a_i,b_i$ given by (\ref{4a1})-(\ref{4bi}) with now $\delta$ at most order $1/\lambda$ singular corresponding to $\tilde k_i,\tilde l_i$ nonsingular. We assume here that $g^{ij}$ and $v^{ij}$ are nonsingular as $\lambda \to 0$ so that $\tilde a_i,\tilde b_i$ are also. The condition in Proposition~\ref{propbetaaibi} for $\delta^2$ to be a left module map gives $v^{ij}$ in terms of $g^{ij}$ as before and the unchanged form
\[ \tilde l_1=-\tilde k_2-\lambda\tilde k_1,\quad  \tilde l_2=\tilde k_1.\]
In what follows will limit ourselves to this case,  where $g^{ij}$ are given, $\tilde k_i$ are our parameters for the connection and everything else is determined. 

\begin{theorem}\label{thmbetalc}
Let $\delta^2$ is a left-module map and $\delta$ have at most an order $\frac{1}{\lambda}$ singularity. Let $\nabla$ be the connection emerging from the extension data for any $g^{ij}$ and parameters $\tilde k_i$.
\begin{enumerate}
\item  The classical limit of the connection exists and has cotorsion and torsion
\begin{align*}C_1^{cl}&=-\frac{1}{|g|}\left(\frac{2}{r}(g^{11}t+g^{12}+g^{21})+\tilde{k_1}\right), \quad C_2^{cl}=\frac{2g^{11}}{|g|r},\\ 
T_1^{cl}&=-\frac{1}{2|g|}g^{11}\left(\frac{2}{r}(g^{11}t+g^{12})+\tilde{k_1}\right),\\
T_2^{cl}&=\frac{1}{2|g|}\left((g^{12}-g^{21})\left(\frac{2}{r}(g^{11}t+g^{12})+\tilde{k_1}\right)-g^{11}\left(\frac{2}{r}(g^{21}t+g^{22})+\tilde{k_2}\right)-\frac{4|g|}{r}\right).\end{align*}
\item The full connection and its torsion can be written in terms of cotorsion as
\begin{align*}
\nabla e_1=&\frac{g^{12}|g|}{(g^{11})^2}\tilde{C_2}e_1\tens e_1-\frac{|g|}{g^{11}}\tilde{C_2}e_2\tens e_1\\
\nabla e_2=&-\frac{|g|}{(g^{11})^2}(g^{12}\tilde{C_1}-g^{22}\tilde{C_2})e_1\tens e_1-\frac{g^{12}}{(g^{11})^2}\tilde{C_2}e_1\tens e_2\\
&+\frac{|g|}{(g^{11})^2}(g^{11}\tilde{C_1}-g^{21}\tilde{C_2})e_2\tens e_1+\frac{|g|}{g^{11}}\tilde{C_2}e_2\tens e_2\\
T_1=&\frac{1}{g^{11}}\tilde C_2,\quad T_2=\frac{|g|}{(g^{11})^2}(g^{11}\tilde C_1+(g^{12}-g^{21}+g^{11})\tilde C_2)-\frac{2}{r}\end{align*} 
where  
\begin{align*}\tilde{C_1}&=\frac{(g^{11})^2}{2|g|}\left(\frac{2}{|g|r}(g^{21}t+g^{22})+\frac{\tilde{k_2}}{|g|}-\frac{2\lambda}{r}(g^{11}t+g^{12})-\lambda\tilde{k_1}\right)\\ 
\tilde{C_2}&=-\frac{(g^{11})^2}{2|g|}\left(\frac{2}{r}(g^{11}t+g^{12})+\tilde{k_1}\right);\quad \tilde{C_i}={C_i-C_i^{cl}\over \lambda}.\end{align*}

\item $C_i=0$ occurs in our moduli space if and only if the $\delta^2$ bimodule map condition displayed in part (5) of Proposition~\ref{propbetaaibi} holds. In this case $T_i=0$ also, giving a one parameter moduli space of weak quantum Levi-Civita connections with parameter $\tilde k_2$. 
\item The connections in (3)  are quantum Levi-Civita if and only if in addition $g^{22}=-\frac{\lambda g^{12}}{2}$, and have the form
\begin{align*}
\nabla e_1=&-\frac{1}{r}e_1\tens e_1,\quad \nabla e_2=\left(\frac{1}{r}\left(t-\frac{\lambda}{2}\right)-\frac{\tilde{k_2}}{2g^{12}}\right)e_1\tens e_1-\frac{g^{12}}{r}g
\end{align*}
\end{enumerate} 
\end{theorem}
\proof
(i) From the conditions (\ref{delta2leftab}) for $\delta^2$ to be a left module map we find \[b_2=\frac{\tilde{a_1}}{\lambda}, \quad b_1=-\frac{\tilde{a_2}}{\lambda}-\tilde{a_1}.\]
We then substitute our expressions for $a_i,b_i$ into the formulae for the full quantum connection found in the proof of Lemma~\ref{lemTC} to get, 
\begin{align*}
\nabla e_1=&-\frac{1}{2|g|}g^{12}\tilde{a_1}e_1\tens e_1+\frac{1}{2|g|}g^{11}\tilde{a_1}e_2\tens e_1\\
\nabla e_2=&-\frac{1}{2|g|}((g^{22}+\lambda g^{12})\tilde{a_1}+g^{12}\tilde{a_2})e_1\tens e_1+\frac{1}{2|g|}g^{12}\tilde{a_1}e_1\tens e_2\\
&+\frac{1}{2|g|}((g^{21}+\lambda g^{11})\tilde{a_1}+g^{11}\tilde{a_2})e_2\tens e_1-\frac{1}{2|g|}g^{11}\tilde{a_1}e_2\tens e_2.
\end{align*}
The braiding map in this case becomes,
\begin{align*}
\sigma_{e_1}(e_1\tens e_1)=&\frac{1}{2\lambda}[r,\tilde{a_1}]e_1\\
\sigma_{e_1}(e_1\tens e_2)=&\frac{1}{2\lambda}([t,\tilde{a_1}]r-[r,\tilde{a_1}](t-\lambda))e_1\\
\sigma_{e_2}(e_1\tens e_1)=&\frac{1}{2\lambda}([r.\tilde{a_1}](1+\lambda)+[r,\tilde{a_2}])e_1+\frac{1}{2\lambda}[r,\tilde{a_1}]e_2\\
\sigma_{e_2}(e_1\tens e_2)=&\frac{1}{2\lambda}(([r,\tilde{a_2}]+\lambda[r,\tilde{a_1}])(\lambda-t)+([t,\tilde{a_2}]+\lambda[t,\tilde{a_1}])r)e_1\\
\end{align*}
\begin{align*}
\sigma_{e_1}(e_2\tens e_1)=&-\frac{1}{2}[r,\tilde{a_1}]e_1+\frac{1}{2\lambda}[r,\tilde{a_1}]e_2\\
\sigma_{e_1}(e_2\tens e_2)=&\frac{1}{2}([r,\tilde{a_1}](t-\lambda)-[t,\tilde{a_1}]r)e_1+\frac{1}{2\lambda}([t,\tilde{a_1}]r-[r,\tilde{a_1}](t-\lambda))e_2\\
\sigma_{e_2}(e_2\tens e_1)=&\frac{1}{2}([r,\tilde{a_2}]+\lambda[r,\tilde{a_1}])e_1+\frac{1}{2\lambda}([r,\tilde{a_2}]-[r,\tilde{a_1}])e_2\\
\sigma_{e_2}(e_2\tens e_2)=&\frac{1}{2}(([r,\tilde{a_2}]+\lambda[r,\tilde{a_1}])(\lambda-t)+([t,\tilde{a_2}]+\lambda[t,\tilde{a_1}])r)e_1\\
&+\frac{1}{2\lambda}([t,\tilde{a_2}]r-[r,\tilde{a_2}](t-\lambda)-\lambda[t,\tilde{a_1}]r+\lambda[r,\tilde{a_1}](t-\lambda))e_2
\end{align*}
The connection is clearly non-singular and has a classical limit given by, 
\begin{align*}
\nabla e_1=&-\frac{1}{2|g|}g^{12}\left(\frac{2}{r}(g^{11}t+g^{12})+\tilde{k_1}\right)e_1\tens e_1+\frac{1}{2|g|}g^{11}\left(\frac{2}{r}(g^{11}t+g^{12})+\tilde{k_1}\right)e_2\tens e_1\\
\nabla e_2=&-\frac{1}{2|g|}\left(g^{22}\left(\frac{2}{r}(g^{11}t+g^{12})+\tilde{k_1}\right)+g^{12}\left(\frac{2}{r}(g^{21}t+g^{22})+\tilde{k_2}\right)\right)e_1\tens e_1\\
&+\frac{1}{2|g|}g^{12}\left(\frac{2}{r}(g^{11}t+g^{12})+\tilde{k_1}\right)e_1\tens e_2\\
&+\frac{1}{2|g|}\left(g^{21}\left(\frac{2}{r}(g^{11}t+g^{12})+\tilde{k_1}\right)+g^{11}\left(\frac{2}{r}(g^{21}t+g^{22})+\tilde{k_2}\right)\right)e_2\tens e_1\\
&-\frac{1}{2|g|}g^{11}\left(\frac{2}{r}(g^{11}t+g^{12})+\tilde{k_1}\right)e_2\tens e_2
\end{align*}
When $\delta^2$ is a left-module map the cotorsion $C_i$, given in Lemma~\ref{lemTC}, become,
\begin{equation}\label{betaCi} C_1=\frac{1}{|g|}\left(\frac{\lambda^2(g^{11})^2}{2|g|}\left(\frac{\tilde{a_2}}{\lambda}+\tilde{a_1}\right)-\tilde{a_1}-\frac{2g^{21}}{r}\right), \quad C_2=\frac{g^{11}}{|g|}\left(\frac{2}{r}-\frac{\lambda g^{11} }{2|g|}\tilde{a_1}\right).
\end{equation}
These have classical limits as stated. We repeat this process for torsion, in which case we have, 
\begin{equation}\label{betaTi} T_1=-\frac{1}{2|g|}g^{11}\tilde{a_1},\quad T_2=\frac{1}{2|g|}\left((g^{12}-g^{21})\tilde{a_1}-g^{11}\tilde{a_2}-\frac{4|g|}{r}\right).
\end{equation}
Which again have classical limits as stated. 

(ii) These expressions (\ref{betaCi}) and (\ref{betaTi}) for $C_i,T_i$ are each invertibly related to $\tilde{a_i}$, in particular 
\begin{equation}\label{tila1C} \tilde{a_1}=\frac{2|g|}{\lambda(g^{11})^2}(C_2^{cl}-C_2)=-\frac{2|g|}{(g^{11})^2}\tilde{C_2}
\end{equation}
\begin{equation}\label{tila2C} \tilde{a_2}=\frac{2|g|^2}{\lambda(g^{11})^2}((C_1-C_1^{cl})+\lambda(C_2-C_2^{cl}))=\frac{2|g|^2}{(g^{11})^2}(\tilde{C_1}+\lambda\tilde{C_2})
\end{equation} 
where $C_i^{cl}$ are the classical values for the cotorsion as given above. We can then substitute (\ref{tila1C}) and (\ref{tila2C}) into the formulae for the full quantum connection to arrive at the form stated. Furthermore, we can use (\ref{tila1C}) and (\ref{tila2C}) in (\ref{betaTi}) to achieve results similar to that in Section~\ref{Secalpha} whereby we obtained a relationship between the torsion and cotorsion as stated. Note that now the cotorsion coefficients here are not constants and have a particular form in terms of our actual parameters, as stated. 

(iii) From (\ref{betaCi}), we can clearly see that $C_2=0$ if and only if $g^{11}=0$. We then have that $$C_1=-\frac{1}{|g|}\left(\tilde{a_1}+\frac{2g^{21}}{r}\right).$$ Therefore, $C_1=0$ if and only if $\tilde{a_1}=-\frac{2g^{21}}{r}$. We can then use equation ({\ref{tildeai}) to expand $\tilde{a_1}$ to arrive at the conditions $g^{11}=0$, $g^{21}+g^{12}=0$, $\tilde k_1=0$ which are precisely the $\delta^2$ bimodule map conditions displayed in part (5) or Proposition~\ref{propbetaaibi}. It is easy to then substitute these conditions on $g^{ij}$ into (\ref{betaTi}) to see that $T_i=0$ also in this case. We can also write 
\[\tilde{a_1}=\frac{2g^{12}}{r}, \quad \tilde{a_2}=\frac{2}{r}(g^{22}-g^{12}t)+\tilde{k_2}\]
in which case our weak quantum Levi-Civita bimodule connections becomes
\begin{align*}
\nabla e_1=&-\frac{1}{r}e_1\tens e_1,\quad \nabla e_2=\left(\frac{1}{r}(t-\lambda)-\frac{g^{22}}{g^{12}r}-\frac{\tilde{k_2}}{2g^{12}}\right)e_1\tens e_1-\frac{g^{12}}{r}g
\end{align*}
\[\sigma(e_1\tens e_1)=e_1\tens e_1,\quad \sigma(e_1\tens e_2)=e_2\tens e_1,\quad \sigma(e_2\tens e_1)=e_1\tens e_2+\lambda e_1\tens e_1\]
\[ \sigma(e_2\tens e_2)= e_2\tens e_2-\lambda\left({2g^{22}\over g^{12}}+\lambda\right)e_1\tens e_1+\lambda e_1\tens e_2-\lambda e_2\tens e_1\]
where the latter are obtained from  $\sigma(e_i\tens e_j)=g^1\tens\sigma_{g^2}(e_i\tens e_j)$ and 
\begin{align*}
\sigma_{e_1}(e_1\tens e_1)=&0,\quad \sigma_{e_1}(e_1\tens e_2)=g^{12}e_1,\quad \sigma_{e_2}(e_1\tens e_1)=-g^{12}e_1,\quad \sigma_{e_2}(e_1\tens e_2)=g^{22}e_1,\\
\sigma_{e_1}(e_2\tens e_1)=&0,\quad \sigma_{e_1}(e_2\tens e_2)=-\lambda g^{12}e_1+ g^{12}e_2,\quad \sigma_{e_2}(e_2\tens e_1)=-\lambda g^{12}e_1-g^{12}e_2,\\
\sigma_{e_2}(e_2\tens e_2)=&\lambda(g^{22}+\lambda g^{12})e_1+(g^{22}-\lambda g^{12})e_2
\end{align*}

(iv) For metric compatibility we must have 
\[\nabla g=(\nabla\tens\id)g+(\sigma\tens\id)(\id\tens\nabla)g=0.\]
where at this point $$g=\frac{1}{|g|}(g^{22}e_1\tens e_1-g^{12}e_1\tens e_2+g^{12}e_2\tens e_1)$$ and that $|g|=g^{12}g^{12}$. Using this we compute for our weak quantum Levi-Civita connections that 
\begin{align*}
(\nabla\tens\id)g=&-\frac{2g^{22}}{g^{12}g^{12}r}e_1\tens e_1\tens e_1+\frac{1}{g^{12}r}e_1\tens e_1\tens e_2+\frac{1}{g^{12}r}e_1\tens e_2\tens e_1-\frac{1}{g^{12}r}e_2\tens e_1\tens e_1\\
&+\frac{1}{g^{12}}\left(\frac{1}{r}(t-\lambda)-\frac{g^{22}}{g^{12}r}+\frac{\lambda k_2}{2g^{12}g^{12}}\right)e_1\tens e_1\tens e_1\\
(\sigma\tens\id)(\id\tens\nabla)g=&-\frac{1}{g^{12}}\left(\frac{1}{r}(t-\lambda)-\frac{g^{22}}{g^{12}r}+\frac{\lambda k_2}{2g^{12}g^{12}}\right)\sigma(e_1\tens e_1)\tens e_1\\
&-\frac{1}{g^{12}r}\sigma(e_1\tens e_1)\tens e_2+\frac{1}{g^{12}r}\sigma(e_1\tens e_2)\tens e_1\\
&-\frac{1}{g^{12}r}\sigma(e_2\tens e_1)\tens e_1.
\end{align*}
Substituting the values of $\sigma$ and combining, we arrive at the requirement 
\[-\frac{1}{|g|r}(2g^{22}+\lambda g^{12})e_1\tens e_1\tens e_1=0\] which gives us the result stated. \endproof

\begin{lemma}\label{lemcurv} For the connection in Theorem~\ref{thmbetalc}, \begin{enumerate}
\item The curvature for general $g^{ij}$ has classical limit 
\begin{align*}
R_{\nabla}^{cl}(e_1)=&-\frac{1}{|g|}g^{11}\left(\frac{2}{r}(g^{11}t+g^{12})+\tilde{k_1}\right)\frac{1}{r}\Vol\tens e_1\\
R_{\nabla}^{cl}(e_2)=&\frac{1}{|g|}\left(\frac{1}{2}\left(\frac{2}{r}(g^{11}t+g^{12})+\tilde{k_1}\right)\left(\frac{2}{r}(g^{11}t+g^{12}+g^{21})+\tilde{k_1}\right)+g^{11}\left(\frac{2}{r}(g^{21}t+g^{22})+\tilde{k_2}\right)\frac{1}{r}\right)\Vol\tens e_1\\
&-\frac{g^{11}}{|g|}\left(\frac{2}{r}(g^{11}t+g^{12})+\tilde{k_1}\right)\frac{1}{r}\Vol\tens e_2
\end{align*}
\item The one parameter moduli space of weak quantum Levi-Civita connections in case (3) of Theorem~\ref{thmbetalc} are all flat.
\end{enumerate}
\end{lemma}

\proof (i) We begin by first computing the full quantum curvature of the connection assuming it is a left module map using the expression for the connection given in Theorem~\ref{thmbetalc} in terms of the residue functions $\tilde{a_i}$. Recall that quantum curvature is given by \[R_{\nabla}=(\extd\tens\id-\id\wedge\nabla)\nabla\]
We then have
\begin{align*}
R_{\nabla}(e_1)=&\frac{1}{2|g|}g^{11}\tilde{a_1}\left(\frac{\lambda}{2|g|}\tilde{a_1}-\frac{2}{r}\right)\Vol\tens e_1-\frac{1}{2|g|}\extd\tilde{a_1}(g^{12}e_1\tens e_1-g^{11}e_2\tens e_1)\\
R_{\nabla}(e_2)=&-\frac{1}{4|g|^2}\left(\lambda(g^{11})^2[\tilde{a_1},\tilde{a_2}]-2|g|\tilde{a_1}^2-4|g|((g^{21}+\lambda g^{11})\tilde{a_1}+g^{11}\tilde{a_2})\frac{1}{r}\right)\Vol\tens e_1\\
&-\frac{g^{11}\tilde{a_1}}{4|g|^2}\left(\lambda g^{11}\tilde{a_1}+\frac{4|g|}{r}\right)\Vol\tens e_2\\
&-\frac{1}{2|g|}((g^{22}+\lambda g^{12})\extd\tilde{a_1}+g^{12}\extd\tilde{a_2})e_1\tens e_1+\frac{1}{2|g|}g^{12}\extd\tilde{a_1}e_1\tens e_2\\
&+\frac{1}{2|g|}((g^{21}+\lambda g^{11})\extd\tilde{a_1}+g^{11}\extd\tilde{a_2})e_2\tens e_1-\frac{1}{2|g|}g^{11}\extd\tilde{a_1}e_2\tens e_2
\end{align*}
Expanding the $\tilde{a_i}$ according to equations (\ref{tildeai}) in terms of the parameters $\tilde{k_i}$ gives us the full quantum curvature of the connection as
\begin{align*}
R_{\nabla}(e_1)=&\left(\frac{1}{2|g|}g^{11}\left(\frac{2}{r}(g^{11}t+g^{12})+\tilde{k_1}\right)\left(\frac{\lambda}{2|g|}g^{11}\left(\frac{2}{r}(g^{11}t+g^{12})+\tilde{k_1}\right)-\frac{2}{r}\right)-\frac{2\lambda}{|g|r^2}(g^{11})^2\right)\Vol\tens e_1\\
R_{\nabla}(e_2)=&\frac{1}{|g|}\left(\frac{1}{2}\left(\frac{2}{r}(g^{11}t+g^{12})+\tilde{k_1}\right)\left(\frac{2}{r}(g^{11}t+g^{12})+\tilde{k_1}+\frac{2}{r}(g^{21}+\lambda g^{11})\right)-\frac{\lambda^2(g^{11})^2}{r^2}\right)\Vol\tens e_1\\
&-\frac{1}{|g|}\left(\frac{\lambda g^{11}}{r^2}(2g^{21}+\lambda g^{11})-g^{11}\left(\frac{2}{r}(g^{21}t+g^{22})+\tilde{k_2}\right)\frac{1}{r}\right)\Vol\tens e_1\\
&\left(\frac{\lambda(g^{11})^2}{|g|r^2}-\frac{g^{11}}{4|g|^2}\left(\frac{2}{r}(g^{11}t+g^{12})+\tilde{k_1}\right)\left(\lambda g^{11}\left(\frac{2}{r}(g^{11}t+g^{12})+\tilde{k_1}\right)+\frac{4|g|}{r}\right)\right)\Vol\tens e_2
\end{align*}
We can then set $\lambda\rightarrow 0$ to get the classical limit stated.

(ii) Using the above formulae for the full quantum curvature in terms of $\tilde{k_i}$, one can clearly see that setting $g^{11}=0$ means that $R_{\nabla}(e_1)=0$. Setting $g^{11}=0$ and $\tilde{k_1}=0$ gives 
\[R_{\nabla}(e_2)=\frac{1}{|g|}\left(\frac{g^{12}}{r}\left(\frac{2}{r}g^{12}+\frac{2}{r}g^{21}\right)\right)\Vol\tens e_1=0\]
given that $g^{12}=-g^{21}$. \endproof

So far we have focussed on the connection on 1-forms. For the connection applied to forms of degree 2 we have the following lemma:
\begin{lemma}\label{nablavolbeta}
For regular $\delta$ and the covariant derivative $\nabla_{e_i}=\frac{1}{2}\llbracket e_i,\;\rrbracket$, we have that $\nabla_{e_i}\Vol=0$ if and only if $\delta^2$ is a left module map.
\end{lemma}
\proof From the cocycle data given in Lemma~\ref{lemTC} we have
\begin{equation*}\nabla_{e_1}\Vol={1\over 2} \llbracket e_1,\Vol\rrbracket={1\over 2}(b_2-a_1)\Vol,\quad 
\nabla_{e_2}\Vol={1\over 2}\llbracket e_2,\Vol\rrbracket=-{1\over 2}\left(a_2+b_1+\lambda b_2-\frac{2b}{r}\right)\Vol 
\end{equation*}
using $\Vol\perp\Vol$ from Lemma~\ref{perpbic}.\endproof

We also have a Hodge Laplacian defined by $\Delta=\delta \extd +\extd \delta$. We also expand $a_i$ according to Proposition~\ref{propbetaaibi} in order to take the exterior derivative and assume $\delta^2$ a left-module. 
\begin{align*}
\Delta(r)=&\delta\extd r=\delta e_1=a_1=\frac{2}{\lambda r}(g^{11}t+g^{12})+k_1\\
\Delta(t)=&\delta\extd t=\delta\left(\frac{t}{r}e_1\right)+\delta\left(\frac{1}{r}e_2\right)\\
=&\extd\left(\frac{t}{r}\right)\perp e_1+\frac{t}{r}\delta e_1+\extd\left(\frac{1}{r}\right)\perp e_2+\frac{1}{r}\delta e_2=-\frac{b}{r^2}+\frac{t}{r}a_1-\frac{b}{r^2}+\frac{1}{r}a_2\\
=&\frac{2}{\lambda r^2}((g^{11}t+g^{12})(t-\lambda)+(g^{21}t+g^{22}))+\frac{tk_1+k_2}{r}\end{align*}\begin{align*}
\Delta(e_1)=&\extd\delta e_1=\extd a_1=\frac{2g^{11}}{\lambda}\extd\left(\frac{t}{r}\right)+\frac{2g^{12}}{\lambda}\extd\left(\frac{1}{r}\right)=\frac{2g^{11}}{\lambda r^2}e_2-\frac{2g^{12}}{\lambda r^2}e_1\\
\Delta(e_2)=&\delta\extd e_2+\extd\delta e_2=\delta\left(\frac{2}{r}\Vol\right)+\extd b\\
=&\extd\left(\frac{2}{r}\right)\perp\Vol+\frac{2}{r}\delta(\Vol)+\frac{2g^{21}}{\lambda}\extd\left(\frac{t}{r}\right)+\frac{2g^{22}}{\lambda}\extd\left(\frac{1}{r}\right)+2A\extd\left(\frac{1}{r}\right)\\
=&-\frac{2}{r^2}e_1\perp\Vol+\frac{2}{r}b_1e_1+\frac{2}{r}b_2e_2+\frac{2g^{21}}{\lambda r^2}e_2-\frac{2g^{22}}{\lambda r^2}e_1-\frac{2b}{r^2}e_1\\
=&-\frac{2}{r}\left(\frac{2}{\lambda r}(g^{21}t+g^{22})+\frac{2}{r}(g^{11}t+g^{12})+\frac{g^{22}}{\lambda r}+k_2+\lambda k_1\right)e_1\\
&+\frac{2}{r}\left(\frac{2}{\lambda r}(g^{11}t+g^{12})+\frac{g^{21}}{\lambda r}+k_1\right)e_2\end{align*}\begin{align*}
\Delta\Vol=&\delta\extd\Vol+\extd\delta\Vol=(\extd b_1)e_1+(\extd b_2)e_2+2b_2\frac{1}{r}\Vol\\
=&-\left(\extd a_2+\lambda\extd a_1-2b\extd\left(\frac{1}{r}\right)\right)e_1+(\extd a_1)e_2+2a_1\frac{1}{r}\Vol\\
=&\left(\frac{2g^{12}}{\lambda r^2}+\frac{2g^{21}}{\lambda r^2}+\frac{4g^{11}t}{\lambda r^2}+\frac{4g^{11}}{r^2}+\frac{2k_1}{r}\right)\Vol
\end{align*} 

Since our $\delta$ does not have a classical limit for generic $g^{ij}$ there is no question that it coincides with the `geometric codifferential'. For completeness, this comes out as
\begin{proposition}
In the classical limit, the geometric codifferential arising from the extension data via the connection is given by
\begin{align*}
\ci\circ\nabla e_1=&\frac{g^{11}(g^{21}-g^{12})}{|g|r}(g^{11}t+g^{12})\\
\ci\circ\nabla e_2=&\frac{1}{|g|r}((g^{21}g^{21}+g^{12}g^{12}-2g^{11}g^{22})(g^{11}t+g^{12})+g^{11}(g^{21}-g^{12})(g^{21}t+g^{22}))\\
\ci\circ\nabla\Vol=&0
\end{align*}
\end{proposition}
\proof In order to attain a unique classical limit we make use of Theorem~\ref{thmbetalc} and therefore assume $\delta^2$ is left module map. We then apply $\ci$ to the resulting classical connection given in the proof of Theorem~\ref{thmbetalc}.
\begin{align*}\ci\circ\nabla e_1=&-\frac{1}{|g|r}g^{12}(g^{11}t+g^{12})g^{11}+\frac{1}{|g|r}g^{11}(g^{11}t+g^{12})g^{21}\\
=&\frac{g^{11}(g^{21}-g^{12})}{|g|r}(g^{11}t+g^{12})\\
\ci\circ\nabla e_2=&-\frac{1}{|g|r}(g^{22}g^{11}(g^{11}t+g^{12})+g^{11}g^{12}(g^{21}t+g^{22}))+\frac{1}{|g|r}g^{12}g^{12}(g^{11}t+g^{12})\\
&+\frac{1}{|g|r}(g^{21}g^{21}(g^{11}t+g^{12})+g^{11}g^{21}(g^{21}t+g^{22}))-\frac{1}{|g|r}g^{11}g^{22}(g^{11}t+g^{12})\\
=&\frac{1}{|g|r}((g^{21}g^{21}+g^{12}g^{12}-2g^{11}g^{22})(g^{11}t+g^{12})+g^{11}(g^{21}-g^{12})(g^{21}t+g^{22}))
\end{align*}
We have previously shown that $\nabla_{e_i}\Vol=0$ when $\delta^2$ is a left module map, but for completeness and to obtain previously unseen formulae we compute here $\nabla\Vol$ in order to compute $\ci\circ\nabla\Vol$. Thus we have 
\begin{align*}
\nabla\Vol=&g^{1}\tens\nabla_{g^{2}}\Vol\\
=&\frac{1}{|g|}g^{22}e_1\tens\nabla e_1\Vol-\frac{1}{|g|}g^{12}e_1\tens\nabla e_2\Vol-\frac{1}{|g|}g^{21}e_2\tens\nabla e_1\Vol+\frac{1}{|g|}g^{11}e_2\tens\nabla e_2\Vol\\
=&\frac{1}{2|g|}\left(g^{22}(b_2-a_1)+g^{12}\left(a_2+b_1+\lambda b_2-\frac{2b}{r}\right)\right)e_1\tens\Vol\\
&-\frac{1}{2|g|}\left(g^{21}(b_2-a_1)+g^{11}\left(a_2+b_1+\lambda b_2-\frac{2b}{r}\right)\right)e_2\tens\Vol
\end{align*} were we have used the formula in the proof of Lemma~\ref{nablavolbeta}. Using the left-module conditions one then has that $\nabla\Vol=0$. \endproof

We now look at some specific examples. Our general analysis was for a constant quantum metrics $g^{ij}$ without assuming quantum symmetry.

\begin{example}\label{betabh}
The unique quantum symmetric real quantum metric for this model is given in \cite{BegMa:cqg} and has the form
\[ (g^{ij})=\begin{pmatrix} {1\over 1+B\lambda^2} & 0 \cr {\lambda\over 1+B\lambda^2} & {1\over B}\end{pmatrix} \]
and ask for $\delta^2$ a left module map, which fixes $v^{ij}$ and $\tilde l$ with
\[ (v^{ij})=\begin{pmatrix} {\frac{-2\lambda}{1+B\lambda^2}} & {\frac{1}{1+B\lambda^2}} \cr {-\frac{1}{B}} & 0\end{pmatrix}, \] leaving
\[\tilde{a_1}=\frac{2t}{r(1+B\lambda^2)}+\tilde{k_1},\quad \tilde{a_2}=\frac{2\lambda t}{r(1+B\lambda^2)}+\frac{2}{rB}+\tilde{k_2}\]
and remaining parameters $\tilde k_i$. We cannot apply parts (3),(4) of the Theorem~\ref{thmbetalc} due to the form of the metric. In fact one has 
\begin{align*}
\delta^2(\Vol)=&-\frac{4t^2}{\lambda r^2(1+B\lambda^2)^2}+\frac{4(1-\lambda Bt)}{\lambda r^2B(1+B\lambda^2)}+\frac{2b(2t+3\lambda)}{\lambda r^2(1+B\lambda^2)}-\frac{4\tilde{k_1}}{\lambda r(1+B\lambda^2)}+\frac{2b\tilde{k_1}}{\lambda r}+\frac{\tilde{k_1}^2}{\lambda}
\end{align*} which is clearly far from being central. This confirms that $\delta^2$ is not a bimodule map. 

The quantum connection arising from $(\delta,\perp)$ is therefore not even weak quantum Levi-Civita.  It is given in terms of cotorsion according to Theorem~\ref{thmbetalc} as
\begin{align*}
\nabla e_1=&-\frac{\tilde{C_2}}{B}e_2\tens e_1\\
\nabla e_2=&\frac{1+B\lambda^2}{B^2}\tilde{C_2}e_1\tens e_1+\frac{1}{B}(\tilde{C_1}-\lambda\tilde{C_2})e_2\tens e_1+\frac{\tilde{C_2}}{B}e_2\tens e_2
\end{align*}
where 
\[\tilde{C_1}=\frac{B}{2(1+B\lambda^2)}\left(B(1+B\lambda^2)\left(\frac{2\lambda t}{r(1+B\lambda^2)}+\frac{2}{rB}+\tilde{k_2}\right)-\frac{2\lambda t}{r(1+B\lambda^2)}-\lambda\tilde{k_1}\right)\]
\[\tilde{C_2}=-\frac{B}{2(1+B\lambda^2)}\left(\frac{2t}{r(1+B\lambda^2)}+\tilde{k_1}\right)\]

The connection has classical limit
\begin{align*}
\nabla e_1=&\frac{B}{2}\left(\frac{2t}{r}+\tilde{k_1}\right)e_2\tens e_1\\
\nabla e_2=&-\frac{1}{2}\left(\frac{2t}{r}+\tilde{k_1}\right)e_1\tens e_1+\frac{B}{2}\left(\frac{2}{rB}+\tilde{k_2}\right)e_2\tens e_1-\frac{B}{2}\left(\frac{2t}{r}+\tilde{k_1}\right)e_2\tens e_2
\end{align*}
where Theorem~\ref{thmbetalc}  gives us the classical torsion and cotorsion as
\begin{align*}
C_1^{cl}=&-B\left(\frac{2t}{r}+\tilde{k_1}\right),\quad C_2^{cl}=\frac{2B}{r}\\
T_1^{cl}=&-\frac{B}{2}\left(\frac{2t}{r}+\tilde{k_1}\right),\quad T_2^{cl}=-\frac{B}{2}\left(\frac{6}{rB}+\tilde{k_2}\right)
\end{align*}
Lemma~\ref{lemcurv} gives the classical limit of the curvature as
\begin{align*}
R_{\nabla}^{cl}(e_1)=&-\frac{B}{r}\left(\frac{2t}{r}+\tilde{k_1}\right)\Vol\tens e_1\\
R_{\nabla}^{cl}(e_2)=&B\left(\frac{1}{2}\left(\frac{2t}{r}+\tilde{k_1}\right)^2+\frac{1}{r}\left(\frac{2}{rB}+\tilde{k_2}\right)\right)\Vol\tens e_1-\frac{B}{r}\left(\frac{2t}{r}+\tilde{k_1}\right)\Vol\tens e_2
\end{align*}
The classical Ricci tensor here is not proportional to the metric (and nor would we expect it to be as the connection is not the Levi-Civita one).

The quantum Laplacian has formulae
\begin{align*}
\Delta(r)=&\frac{2t}{\lambda r(1+B\lambda^2)}+k_1\\
\Delta(t)=&\frac{2}{\lambda r^2}\left(\frac{t^2}{1+B\lambda^2}+\frac{1}{B}\right)+\frac{tk_1+k_2}{r}\\
\Delta(e_1)=&\frac{2}{\lambda r^2(1+B\lambda^2)}e_2\\
\Delta(e_2)=&\left(-\frac{8t}{r^2(1+B\lambda^2)}-\frac{6}{\lambda r^2B}-\frac{2(\lambda k_1+k_2)}{r}\right) e_1\\
&+\left(\frac{4t}{\lambda r^2(1+B\lambda^2)}+\frac{2k_1}{r}+\frac{2}{r^2(1+B\lambda^2)}\right)e_2
\end{align*}
and like $\delta$ does not have a classical limit.
\end{example}

\begin{example}\label{betaglc} Clearly the nicest form of the metric in the sense that all the cases of Theorem~\ref{thmbetalc} hold, is 
\[ (g^{ij})={1\over B}\begin{pmatrix} 0 & 1 \cr -1 & {-\frac{\lambda }{2}} \end{pmatrix}\]
for an  overall normalisation $B$, and we also assume that $\tilde k_1=0$ for part (3) of the theorem to apply and $C_i=T_i=0$. From the formulae displayed in (\ref{tildeai}) we have 
\[\tilde{a_1}=\frac{2}{rB},\quad \tilde{a_2}=-\frac{2t}{rB}-\frac{\lambda}{rB}+\tilde{k_2}.\] 
From Proposition~\ref{propbetaaibi}, for $\delta^2$ a left module map we have that 
\[ (v^{ij})=\frac{1}{B}\begin{pmatrix} {1} & 0 \cr {-\frac{\lambda}{2}} & {1}\end{pmatrix}, \] 
and also fix $\tilde b_i$. From this data one can compute $\delta^2(\Vol)=0$ so that $\delta^2$ is a bimodule map, as it must be according to Proposition~\ref{propbetaaibi} for this form of metric. The 1-parameter family of quantum Levi-Civita connections according to Theorem~\ref{thmbetalc} are then given by 
\begin{align*}
\nabla e_1=&-\frac{1}{r}e_1\tens e_1,\quad \nabla e_2=\left(\frac{t}{r}-\frac{B\tilde{k_2}}{2}\right)e_1\tens e_1+\frac{1}{r}e_1\tens e_2-\frac{1}{r}e_2\tens e_1
\end{align*}
with braiding map 
\[\sigma(e_1\tens e_1)=e_1\tens e_1,\quad \sigma(e_1\tens e_2)=e_2\tens e_1,\quad \sigma(e_2\tens e_1)=e_1\tens e_2+\lambda e_1\tens e_1\]
\[ \sigma(e_2\tens e_2)=e_2\tens e_2+\lambda(e_1\tens e_2-e_2\tens e_1)\]
and zero curvature by Lemma~\ref{lemcurv}. 

Finally, we have the Hodge Laplacian given by
\begin{align*}
\Delta(r)=&0,\quad \Delta(t)=-\frac{3}{r^2B}+\frac{k_2}{r},\quad \Delta(e_1)=-\frac{2}{\lambda r^2B}e_1\\
\Delta(e_2)=&-\frac{2}{r}\left(\frac{1}{2rB}-\frac{2t}{\lambda rB}+k_1\right)e_1-\frac{2}{r}\left(\frac{1}{\lambda rB}+k_1\right)e_2,\quad \Delta(\Vol)=0.
\end{align*} 
\end{example}

\begin{example}\label{betainner} Since the calculus is inner with $\theta=-\frac{\extd t}{\lambda}$, we also have a canonical example of $\delta$ with $\cj=\frac{1}{2}\perp$, in particular $\cj_{e_i}(e_j)=g^{ij}$ is similar to the preceding example but now with 
\[  (g^{ij})=\frac{b}{2}\begin{pmatrix} 0 & 1 \cr -1 & -\lambda\end{pmatrix}. \]
We also compute $\delta(e_i)=\theta\perp e_i$ and $\delta\Vol=\theta\perp\Vol$  as 
\[\delta(e_1)=\frac{b}{\lambda r},\quad \delta(e_2)=\frac{b}{r}(1-\frac{t}{\lambda }),\quad \delta(\Vol)=\theta\perp\Vol=\frac{b}{\lambda r}(te_1+e_2)\]
This corresponds to parameters  $\tilde k_i=0$ and 
\[\tilde{a_1}=\frac{b}{r}, \quad \tilde{a_2}=-\frac{b}{r}(t+\lambda ),\quad \tilde b_1=\frac{b}{r}t=-\tilde{a_2}-\lambda\tilde{a_1},\quad \tilde b_2=\frac{b}{r}=\tilde{a_1}\]
From this or from  $\cj_{e_i}\Vol={1\over 2}e_i\perp\Vol$ to compute the form of $v^{ij}$ we see that $\tilde l_i=0$ so that $\delta^2$ is a left module map by our analysis. Furthermore one can check that $\delta^2(\Vol)=0$ so that $\delta^2$ is a bimodule map as it must be according to Proposition~\ref{propbetaaibi}. Formulae in the proof of Theorem~\ref{thmbetalc} allow us to compute the connection from $\tilde a_i$, as 
\begin{align*}
\nabla e_1=&-\frac{1}{r}e_1\tens e_1,\quad \nabla e_2=\frac{1}{r}((t+\lambda)e_1\tens e_1+e_1\tens e_2-e_2\tens e_1)
\end{align*}
\[\sigma(e_1\tens e_1)=e_1\tens e_1,\quad \sigma(e_1\tens e_2)=e_2\tens e_1,\quad \sigma(e_2\tens e_1)=e_1\tens e_2+\lambda e_1\tens e_1\]
\[ \sigma(e_2\tens e_2)=e_2\tens e_2+\lambda(e_1\tens e_2-e_2\tens e_1+\lambda e_1\tens e_1)\]
and Theorem~\ref{thmbetalc} tells is that this is torsion free and cotorsion free or `weak quantum Levi-Civita'. It is flat but not fully quantum Levi-Civita since $g^{22}\neq-\frac{\lambda}{2}g^{12}$, in fact 
\[\nabla g=\frac{2\lambda}{br}e_1\tens e_1\tens e_1\]
so that the classical limit is metric compatible. 

The quantum Laplacian, given by $\Delta=\delta\extd+\extd\delta$, is
\begin{align*}
\Delta(r)=&\frac{b}{\lambda r},\quad \Delta(t)=0, \\ \Delta(e_1)=&-\frac{b}{\lambda r^2}e_1,\quad
\Delta(e_2)=\left(\frac{2bt}{\lambda r^2}+\frac{b}{r^2}\right)e_1+\frac{b}{\lambda r^2}e_2, \quad \Delta(\Vol)=0.
\end{align*}
The quantum Laplacian here is singular so does not have a classical limit, as for the codifferential.  \end{example}

\section{Conclusions and discussion}\label{Secsemi}

We have seen that the new approach to classical Riemannian geometry and its quantisation in \cite{Ma:rec} via an axiomatic `codifferential'
$\delta$ works very well for the $\alpha$ calculus on our quantum spacetime (Section~3) and does give the quantum Levi-Civita connection for this model when $g$ is quantum symmetric as assumed in \cite{MaTao}. One may expect that this will also be the case
for  other quantum differential spacetimes that are in some (to be determined) sense `close enough' to classical. 

It is also striking that in both cases asking for $\delta^2$ to be a left module map or `left-tensorial', in the sense $\delta^2(f\Vol)=f\delta^2\Vol$ for all $f$ in the quantum spacetime algebra, ensures that
the connection coming from $(\delta,\perp)$ in our quantum Koszul formula is nonsingular as $\lambda\to 0$ (more generally, it needs to hold at least to order $\lambda$). We also saw how this left module map property links the induced interior product $\cj$ to the metric extended as something like a derivation, possibly with $O(\lambda)$ corrections. And we saw that in both cases $\perp$ does not have to be symmetric even though that would be the classical choice (where $\perp=(\ ,\ )$ (the metric) extended in both arguments to forms). In Section
3 we saw that the symmetric choice allows $\delta$ to agree with the geometric divergence defined as $(\ ,\ )\nabla$ while in Section~4
only an antisymmetric plus $O(\lambda)$ choice was allowed by the differential calculus, which is a first hint that it is in some sense `far from classical'. Finally, we saw in both cases how $\delta^2$ being additionally  a right module map or  `right-tensorial' (hence a bimodule map) is a further constraint which in Section~3 forces the metric to be symmetric and lands us on the quantum Levi-Civita connection, while in Section~4 it forces the metric to be mostly antisymmetric (leaving $g^{22}$ unconstrained) and lands us on a weak quantum Levi-Civita connection as in  Example~\ref{betaglc}. Requiring this to be fully quantum Levi-Civitia then fixes the relative value of $g^{22}$ also. Thus we are forced to a form of metric that is not symmetric but {\em antisymmetric} in the classical limit. In other words, the quantum Koszul formula method which we have explored works also for the $\beta$ calculus model on our quantum spacetime in Section 4 but the geometry that it quantises more naturally is symplectic rather than Riemannian. 
It is fair to say that this huge contrast was not visible until now, where both models have been studied as different quantum Riemannian geometries of not fundamentally different character if one just wants a quantum symmetric metric and quantum Levi-Civita connection\cite{MaTao,BegMa:cqg}. The difference now is that we want the geometry to emerge as part of a quantisation of connections and interior products on higher differential forms as well as on $\Omega^1$, which is an integral part of the the quantum Koszul formula, i.e. we want the quantum-`Riemannian' geometry to work with differential forms in the spirit of Hodge theory and the Cartan formula for codifferentials.

 It is not clear of course if our in-depth analysis of one particular spacetime $[r,t]=\lambda r$ allows us to draw lessons more widely. The above phenomena would need to explored in other models; suffice it to say that some of these general features echo some of the steps in proof in \cite{Ma:rec}  that we can recover classical Riemannian geometry from axiomatic properties of $\delta$ of classical type. It should also be pointed out that the central extension formalism in \cite{Ma:rec} of which the 
$(\delta,\perp)$ construction is an example is more general and there could be other constructions leading to flat central extensions. 
Moreover, it seems likely that the central extension theory should itself be generalised in order to recover the actual $\beta=1$ quantum Riemannian geometry in \cite{BegMa:cqg}. This is because the differential calculus on this model has in fact a natural one higher-dimension extension dictated by quantum Poincar\'e group invariance\cite{Sit}. Namely in 2D this is the 3D calculus  with
\[  [\extd r,r]=\lambda\theta',\quad [\extd r,t]=0,\quad[\extd t,t]=\lambda(\theta'-\extd t),\quad [\theta',r]=0,\quad[\theta',t]=\lambda\theta' \]
which we see is {\em not} a central extension.  Rather, it is shown in \cite{Ma:rec} that this calculus is more like a central extension of the calculus on $r$ followed by a semidirect product construction along the lines \cite{Ma:alm}. This in turn works more generally for quantum spacetimes of the form $C^\infty(N)\lcross \R$ where $N$ is a spatial Riemannian manifold and the semidirect product of space with a time coordinate is given by the action of a conformal killing vector. Such quantum spacetimes we called `almost commutative' in \cite{Ma:alm} and the $\beta=1$ calculus is an example in this family with Killing vector $r{\del\over\del r}$. Therefore a direction for further work could be to extend the analysis of Section~4 to 
the quantum Koszul construction for this more general class. It would also be interesting to explore it for finite groups where several quantum Riemannian geometries in our sense are known, as well as for $q$-deformed examples such as $q$-$SU_2$ and the $q$-sphere. 

Finally, one should continue the process of making contact between constructive approaches and  other more `top down' (but more powerful) approaches to noncommutative geometry, most notably that of Connes \cite{Con} based on an axiomatic Dirac operator $D$ or `spectral triple' rather than $\delta$. One could perhaps consider $\extd+\delta$ in this vein as a first step. Better, one should extend the central extension point of view of \cite{Ma:rec} to include spinors and make proper contact with the actual geometric Dirac operator and its interaction with $\delta$. It would also be interesting to make contact with more recent work such as \cite{Rie}. These are some directions for further work.

\end{document}